\documentclass[reqno,11pt]{amsart}
\usepackage[english]{babel}
\usepackage[utf8]{inputenc}
\usepackage[T1]{fontenc}
\usepackage{lmodern}
\usepackage{graphicx}
\usepackage{amsmath,amssymb,amsthm,bm,mathrsfs}
\usepackage{enumerate,array,tabularx,booktabs,ragged2e,float,adjustbox}
\usepackage{fullpage}
\usepackage[numbers]{natbib}
\usepackage{xcolor}
\definecolor{violet}{RGB}{128,0,160}
\usepackage{hyperref}
\usepackage{todonotes}
\hypersetup{colorlinks=true,
linkcolor=blue,
citecolor=magenta,
urlcolor=blue,
breaklinks=true,
pdfstartview=FitH,
bookmarksopen=true,
pdftitle={Positivity and Boundedness of Solutions to the Guyon--Lekeufack Volatility Model with General Kernels},
pdfauthor={Julien Guyon, Thibault Jeannin, Benjamin Jourdain}}

\newcommand{\E}{\mathbb E}

\makeatletter
\newcommand{\myitem}[2][]{\item[#1]\@hyper@itemfalse\phantomsection\def\@currentlabel{#1}\label{#2}}
\newcommand{\myref}[1]{\hyperref[#1]{\ref*{#1}}}
\makeatother
\numberwithin{equation}{section}
\theoremstyle{plain}
\newtheorem{prooff}{Proof}[section]
\newtheorem{lemma}[prooff]{Lemma}
\newtheorem{cor2}{Corollary}
\newtheorem{prop2}{Proposition}
\newtheorem{thm2}{Theorem}
\theoremstyle{definition}
\newtheorem{rk}[prooff]{Remark}
\newtheorem{example}[prooff]{Example}
\newtheorem{def2}[prooff]{Definition}
\newtheorem{notation}{Notation}
\newtheorem{assumption}{Assumption}
\allowdisplaybreaks

\title{Positivity and
Boundedness of Solutions to the
Guyon--Lekeufack Volatility Model with General Kernels}
\author{Julien Guyon$^{1,2}$, Thibault Jeannin$^1$, Benjamin Jourdain$^1$}
\date{}
\begin{document}
\begin{abstract}
We study nonnegativity of the volatility in the continuous-time Guyon--Lekeufack model, in which volatility is an affine function of the convolution of past returns with a kernel $K_1$ and the square root of the integral of past squared returns weighted by a kernel $K_2$.
We show that, when volatility is nonnegative and $K_1$ admits a resolvent of the first kind, uniform boundedness of integrals of $K_2$ with respect to this resolvent yields a deterministic upper bound for the asset price, which is then a bounded martingale on every finite time interval.
This result is the counterpart of the strict local martingality criterion for negative volatility in \cite[Corollary~2.6(ii)]{GuyonJeanninJourdain2026Martingale}.
We also derive such bounds for popular stochastic volatility models when the Brownian motions driving the asset price and volatility have correlation $-1$.
We give sufficient conditions for nonnegative volatility with a
completely monotone return convolution kernel $K_1$, and necessary conditions when $K_1$ and $K_2$ are respectively $M$- and $N$-exponential kernels.
For $M=1$, we obtain a necessary and sufficient condition.
\end{abstract}

\maketitle
\footnotetext[1]{\raggedright CERMICS, CNRS, ENPC, Institut Polytechnique de Paris,
Marne-la-Vallée, France. Emails:
\href{mailto:julien.guyon@enpc.fr}{julien.guyon@enpc.fr},
\href{mailto:thibault.jeannin@enpc.fr}{thibault.jeannin@enpc.fr},
\href{mailto:benjamin.jourdain@enpc.fr}{benjamin.jourdain@enpc.fr}.}
\footnotetext[2]{\raggedright NYU Tandon School of Engineering, Department of Finance
and Risk Engineering, One MetroTech Center, Brooklyn, NY 11201, USA.
Email: \href{mailto:julien.guyon@nyu.edu}{julien.guyon@nyu.edu}.
This research benefited from the support of the chair \textit{Risques
Financiers}, Fondation du Risque, and from the BNP Paribas chair
\textit{Futures of Quantitative Finance}.}

\bigskip
\noindent\textbf{Acknowledgements:}
The authors thank Paul Invernizzi for insightful discussions.
To sharpen an earlier sufficient condition for nonnegative volatility, the authors proposed using the resolvent of the first kind and a semimartingale decomposition. They used OpenAI's GPT-6 Astra model to help adapt the proof to this approach; the sharp bounds for exponential kernels had already been derived by the authors.

\section{Introduction}
Empirical evidence suggests that the volatility of financial markets depends on the past evolution
of the asset price
\cite{Zumbach2009,Zumbach2010,ChicheporticheBouchaud2014,
Bouchaud2021,AndreBoumezoudJourdain2023}.
Guyon and Lekeufack \cite{GuyonLekeufack2023} introduced a
path-dependent volatility (PDV) model based on past returns and past squared returns
weighted by memory kernels. We study when its volatility\footnote{In this paper, by ``volatility'' we mean $\sigma_t$, not $|\sigma_t|$, in \eqref{PDV} below.} stays nonnegative.
From a modeling perspective, verifying the nonnegativity of volatilities is important because it
ensures that the correlations specified between
the driving Brownian motions retain their interpretation as
instantaneous correlations between asset returns (in multi-asset models) or between asset returns and changes in some extra stochastic processes, such as stochastic interest rates or stochastic variance (in single- or multi-asset models).
For a Brownian motion $W$ and deterministic continuous functions
$g_1:\mathbb R_+\to\mathbb R$ and $g_2:\mathbb R_+\to[0,\infty)$,
we define the PDV model through the coupled system, for $t\ge0$,
\begin{align}
\frac{dS_t}{S_t} &= \sigma_t\, dW_t, \tag{PDV} \label{PDV}\\
\sigma_t &\mathrel{:=} \beta_0+\beta_1 \bar R_{1,t}+\beta_2 \sqrt{\bar R_{2,t}}, \notag\\
\bar R_{1,t} &\mathrel{:=} g_1(t)
        + \int_{0}^{t} K_1(t-s)\sigma_s\, dW_s, \notag\\
\bar R_{2,t} &\mathrel{:=} g_2(t)+\int_0^t K_2(s,t)\sigma_s^{2}\,ds, \notag
\end{align}
where $\beta_0,\beta_2>0$ and $\beta_1<0$.
The nonnegative convolution kernel $K_1:(0,\infty)\to[0,\infty)$
and the nonnegative kernel $K_2:\overline\Gamma\to\mathbb R_+$ weight
past returns and past squared returns, respectively, with
$\overline\Gamma:=\{(s,t)\in\mathbb R_+^2:s\le t\}$.

When $K_1$ and $K_2$ are respectively $M$- and $N$-exponential
kernels (defined below), with integers $M,N\ge1$, the PDV model admits a
finite-dimensional Markovian representation.
The kernels $K_1$ and $K_2$ of the four-factor PDV model
are both $2$-exponential kernels.
Existence and pathwise uniqueness of strong solutions were established for the four-factor model
in \cite[Theorem~1.2]{NutzRiveros2024} and for the PDV model with more general
Volterra kernels in \cite[Theorem~1]{AndresJourdain2023}. For
$K_1(t)=\lambda_{1,1}e^{-\lambda_{1,1}t}$ and
$K_2(s,t)=\lambda_{2,1}e^{-\lambda_{2,1}(t-s)}$,
Guyon and Lekeufack \cite[Section~4.1.5]{GuyonLekeufack2023} observed that
$\lambda_{2,1}<2\lambda_{1,1}$ ensures nonnegativity when
$\sigma_0\ge0$. Nutz and Riveros Valdevenito
\cite[Theorem~2.6]{NutzRiveros2024} obtained a strictly positive pathwise
lower bound when $\sigma_0>0$.
Andr\`es and Jourdain \cite[Theorem~2]{AndresJourdain2023} extended
the positivity result to possibly nonconvolution return kernels
$K_1(s,t)=g(s)e^{h(t)}$ and general squared-return kernels $K_2(s,t)$
satisfying $\partial_tK_2(s,t)-2h'(t)K_2(s,t)\ge0$.

Nonnegativity and invariance for stochastic Volterra equations
are studied in
\cite{Alfonsi2023NonnegativityKernels,AlfonsiSzulda2024,
AbiJaberCuchieroPelizzariPulidoSvalutoFerro2024PolynomialVolterraPolynomialVolterra}.
Bondi and Pulido \cite{BondiPulido2024} establish Feller-type
tests for scalar convolution equations with nonsingular kernels.
Their tests do not directly apply to \eqref{PDV} because
$\sigma$ depends on both $\bar R_1$ and $\bar R_2$, and the initial inputs
depend on time. As in \cite[Lemma~1]{BondiPulido2024},
we derive a semimartingale decomposition for $\sigma$
and obtain sufficient conditions for positivity
by controlling the drift term arising from memory.

Strikingly, for every $T>0$, when $\sigma\ge0$ and condition~\eqref{eq:C2T} below holds, Theorem~\ref{thm:GL_positive_moments_compact} below gives a deterministic upper bound on $\sup_{t\in[0,T]}S_t$, despite the Brownian driving noise. This bound depends only on the kernels, coefficients and initial inputs.
Hence $S$ is a bounded martingale on $[0,T]$.
Beyond martingality, this bound restricts the support of the price
distribution, ensures finiteness of all positive moments, and implies
zero call prices for strikes above the bound.
We also extend this price-bound principle to several stochastic volatility models in the limiting case $\rho=-1$ of perfect negative correlation between the Brownian motions driving the asset price and volatility. In the PDV model, this result and the strict local martingality
criterion in \cite[Corollary~2.6~(ii)]{GuyonJeanninJourdain2026Martingale}
yield a dichotomy according to the sign of volatility:
for smooth return kernels $K_1$ admitting a resolvent of the first kind, $S$ is a strict local martingale on $[0,T]$
whenever $\sigma_{t_0}<0$ has positive probability
for some $t_0\in[0,T)$.

For a completely monotone return kernel $K_1$ (defined below),
Theorem~\ref{prop:pos-MR1-generalK2-vector0} below gives explicit
conditions on the kernels and initial state that ensure nonnegative
volatility, together with an explicit process that bounds $\sigma$
from below. In particular, its specialization to the four-factor PDV
model provides sufficient conditions on both the kernel parameters
and the initial factors, see the second row of
Table~\ref{tab:suff-necessary-conditions} below. It thus addresses the question of finding
sufficient conditions for nonnegative volatility in that model,
left open by Nutz and Riveros Valdevenito
\cite[Section~3.2]{NutzRiveros2024}.

When $K_1$ and $K_2$ are respectively $M$- and $N$-exponential
kernels, Theorem~\ref{gonegative-N-eff2} below gives sufficient
conditions under which volatility is negative with positive
probability at every positive time. By contraposition, these yield
necessary conditions for nonnegative volatility.
For $M=1$, Corollary~\ref{cor:M1-equivalence} below shows that when the model is initialized from past market returns
with $\sigma_0>0$, nonnegative volatility
and martingality of $S$ on every finite time interval
are equivalent to the same explicit condition on the kernel rates.
For $M\ge2$, Table~\ref{tab:suff-necessary-conditions} below displays
the remaining gap between the sufficient and necessary conditions.

\medskip

\section{Setting and well-posedness}
The filtered probability space
$(\Omega,\mathcal F,(\mathcal F_t)_{t\ge0},\mathbb P)$ satisfies the usual conditions,
and $W$ is an $(\mathcal F_t)$-Brownian motion.
Throughout the paper, the initial asset price $S_0>0$ is deterministic.
For the PDV model, the initial values satisfy
\begin{equation}\label{eq:initial-values}
\bar R_{1,0}=g_1(0),\qquad \bar R_{2,0}=g_2(0)>0,\qquad
\sigma_0=\beta_0+\beta_1g_1(0)+\beta_2\sqrt{g_2(0)}.
\end{equation}
\begin{notation}
If $X$ has almost surely locally absolutely continuous paths, $\dot X$ denotes their ordinary time derivative, defined for almost every time on each such path. In particular, almost surely, $X_t=X_0+\int_0^t\dot X_s\,ds$ for every $t\ge0$. For a locally integrable function $f$ and a locally finite measure $\mu$ on $\mathbb R_+$, we define the convolution by $(f*\mu)(t)\mathrel{:=}\int_{[0,t]}f(t-s)\,\mu(ds)$. The same notation applies to the convolution of two functions.
\end{notation}

\subsection{Convolution Kernels and initial inputs}
\label{sec:semimartingale-decomposition}

\begin{def2}
A convolution kernel $K$ on $(0,\infty)$ with values in $[0,\infty)$ is completely monotone if
$K\in C^\infty((0,\infty))$ and
\[
\forall n\in\mathbb N\cup\{0\},\quad \forall t>0,\qquad
(-1)^nK^{(n)}(t)\ge0.
\]
By Bernstein's theorem, complete monotonicity together with
$\int_0^\infty K(t)\,dt=1$ is equivalent to the representation
\[
\forall t>0,\qquad
K(t)=\int_{\lambda\in(0,\infty)}\lambda e^{-\lambda t}\,\nu(d\lambda),
\]
where $\nu$ is a unique probability measure on $(0,\infty)$.
\end{def2}
\begin{assumption}\label{assumption4}
The return kernel $K_1$ is completely monotone and admits the normalized representation
\[
\forall u>0,\qquad
K_1(u)=\int_{\lambda\in(0,\infty)}\lambda e^{-\lambda u}\,\nu(d\lambda),
\]
where $\nu$ is a probability measure on $(0,\infty)$.
Moreover, the second moment of $\nu$ is finite:
\begin{equation}\label{eq:nu-second-moment}
\int_{\lambda\in(0,\infty)}\lambda^2\,\nu(d\lambda)<\infty.
\end{equation}
The initial input $g_1$ has the representation
\begin{equation}\label{eq:initial-spectral-representation}
\forall t\ge0,\qquad
g_1(t)=\int_{\lambda\in(0,\infty)}
e^{-\lambda t}R_{1,0}(\lambda)\,\nu(d\lambda),
\end{equation}
where the deterministic measurable initial factor profile
$R_{1,0}:(0,\infty)\to\mathbb R$ satisfies
\begin{equation}\label{eq:initial-spectral-regularity}
\int_{\lambda\in(0,\infty)}
(1+\lambda)|R_{1,0}(\lambda)|\,\nu(d\lambda)<\infty.
\end{equation}
\end{assumption}
The return component $\bar R_{1,t}$ then admits the factor representation
\[
\begin{aligned}
\forall t\ge0,\qquad
&\bar R_{1,t}=\int_{\lambda\in(0,\infty)}
R_{1,t}(\lambda)\,\nu(d\lambda),\\
\forall t\ge0,\quad\forall\lambda>0,\qquad
&dR_{1,t}(\lambda)=-\lambda R_{1,t}(\lambda)\,dt
+\lambda\sigma_t\,dW_t.
\end{aligned}
\]
Thus, once the kernels and coefficients are fixed, the model
is specified by $S_0$, $g_2$ and the initial factor profile $R_{1,0}(\cdot)$.

\begin{rk}
The normalization $\int_0^\infty K_1(u)\,du=1$ does not restrict the PDV dynamics. Indeed, if $c\mathrel{:=}\int_0^\infty K_1(u)\,du\in(0,\infty)$, replacing $K_1$, $g_1$ and $\beta_1$ by $\frac{K_1}{c}$, $\frac{g_1}{c}$ and $c\beta_1$, respectively, changes $\bar R_1$ into $\frac{\bar R_1}{c}$ and leaves $\sigma$ and $S$ unchanged.
\end{rk}

\begin{notation}\label{not:return-kernel}
Under Assumption~\ref{assumption4}, we set
\[
\Lambda_1\mathrel{:=}\int_{\lambda\in(0,\infty)}\lambda\,\nu(d\lambda)>0,\qquad
\lambda_{1,\mathrm{eff}}
\mathrel{:=}\frac{\int_{\lambda\in(0,\infty)}\lambda^2\,\nu(d\lambda)}{\Lambda_1}.
\]
Condition~\eqref{eq:nu-second-moment} extends $K_1$ canonically to a function in $C^1(\mathbb R_+)$; we use the same symbol for this extension. It also gives $\Lambda_1,\lambda_{1,\mathrm{eff}}<\infty$ and
\[
K_1(0)=\Lambda_1,\qquad K_1'(0)=-\Lambda_1\lambda_{1,\mathrm{eff}}.
\]
\end{notation}

\begin{def2}
For a locally integrable convolution kernel $K$ taking values in $\mathbb R$, a nonnegative locally
finite measure $\mathcal R$ is a resolvent of the first kind if
\[
\text{for almost every }t>0,\qquad
\int_{[0,t]}K(t-s)\,\mathcal R(ds)=1.
\]
\end{def2}
For a nonzero completely monotone kernel $K\in C^1(\mathbb R_+)$,
there exists a unique resolvent of the first kind. It is the sum of $\frac{1}{K(0)}\delta_0$ and an absolutely continuous
measure with a completely monotone density
\cite[Chapter~5, Theorems~5.2 and~5.4]{GripenbergLondenStaffans1990}.
Under Assumption~\ref{assumption4}, the kernel $K_1$ is completely
monotone. We denote by $\mathcal R$ the resolvent of the first kind of $K_1$.
\begin{def2}\label{NexponentialKernel}
An $N$-exponential kernel has the form
\[
\forall u\ge0,\qquad
K(u)=\sum_{j=1}^N\theta_j\lambda_j e^{-\lambda_j u},
\qquad \theta_j>0,\quad \sum_{j=1}^N\theta_j=1,\quad
\lambda_1>\cdots>\lambda_N>0.
\]
It is therefore completely monotone, with
\[
\nu=\sum_{j=1}^N\theta_j\delta_{\lambda_j}.
\]
When $K_1$ is an $M$-exponential kernel and $K_2$ is an $N$-exponential kernel, their weights and rates are denoted by $(\theta_{1,i},\lambda_{1,i})_{1\le i\le M}$ and $(\theta_{2,j},\lambda_{2,j})_{1\le j\le N}$.
The deterministic initial factors satisfy $R_{1,0}(\lambda_{1,i})\in\mathbb R$ for every $i\in\{1,\ldots,M\}$ and $R_{2,0}(\lambda_{2,j})\ge0$ for every $j\in\{1,\ldots,N\}$. The initial inputs are specified by
\begin{equation}\label{eq:exponential-inputs}
\forall t\ge0,\qquad
g_1(t)=\sum_{i=1}^M\theta_{1,i}e^{-\lambda_{1,i}t}R_{1,0}(\lambda_{1,i}),
\qquad
g_2(t)=\sum_{j=1}^N\theta_{2,j}e^{-\lambda_{2,j}t}R_{2,0}(\lambda_{2,j}).
\end{equation}
\end{def2}
In this exponential-kernel setting, $K_2(s,t)=K_2(0,t-s)$ for $0\le s\le t$.
Since $R_{2,0}(\lambda_{2,j})\ge0$ and $g_2(0)=\bar R_{2,0}>0$ by~\eqref{eq:initial-values},
\begin{equation}\label{eq:exponential-input-positivity}
\forall t\ge0,\qquad
g_2(t)\ge e^{-\lambda_{2,1}t}g_2(0)=e^{-\lambda_{2,1}t}\bar R_{2,0}>0.
\end{equation}

\begin{rk}
For an $N$-exponential kernel $K$, its resolvent of the first
kind is given by
\begin{equation}\label{eq:N-exponential-resolvent}
\mathcal R(dt)=\frac{1}{\sum_{j=1}^N\theta_j\lambda_j}\delta_0(dt)
+\left(1+\sum_{p=1}^{N-1}
\frac{e^{-\alpha_p t}}{\displaystyle\alpha_p\sum_{j=1}^N
\frac{\theta_j\lambda_j}{(\lambda_j-\alpha_p)^2}}
\right)dt,
\end{equation}
where, for each $p\in\{1,\ldots,N-1\}$,
$\alpha_p\in(\lambda_{p+1},\lambda_p)$ is the unique solution of
\[
\sum_{j=1}^N\frac{\theta_j\lambda_j}{\lambda_j-\alpha_p}=0.
\]
Indeed, applying \cite[Section~4]{DarrowStepaniants2026} to the Laplace transform of $K$, $\widehat K(z)=\sum_{j=1}^N\frac{\theta_j\lambda_j}{z+\lambda_j}$,
and using $K*\mathcal R=1$ gives, for $z>0$ (so that $z\widehat K(z)>0$),
\[
\widehat{\mathcal R}(z)=\frac{1}{z\widehat K(z)}=\frac{1}{\sum_{j=1}^N\theta_j\lambda_j}+\frac{1}{z}+\sum_{p=1}^{N-1}\frac{1}{\alpha_p\sum_{j=1}^N\frac{\theta_j\lambda_j}{(\lambda_j-\alpha_p)^2}}\frac{1}{z+\alpha_p}.
\]
Taking inverse Laplace transforms yields~\eqref{eq:N-exponential-resolvent}.
\end{rk}

\subsection{Well-posedness}
\begin{assumption}\label{assumption3}
The initial input $g_2$ is strictly positive and locally absolutely continuous.
The squared-return kernel $K_2$ is nonnegative and measurable.
For every $s\ge0$, the map $t\mapsto K_2(s,t)$ is
locally absolutely continuous on $[s,\infty)$.
For every $T>0$, the following bounds hold:
\begin{equation}\label{eq:K2-local-bounds}
\sup_{t\in[0,T]}\sup_{s\in[0,t]}K_2(s,t)<\infty,\qquad
\operatorname*{ess\,sup}_{t\in[0,T]}
\int_0^t|\partial_tK_2(s,t)|\,ds<\infty.
\end{equation}
\end{assumption}

\begin{prop2}\label{prop:existence-uniqueness-moments}
Under Assumptions~\ref{assumption4} and~\ref{assumption3}, $g_1$ belongs to $C^1(\mathbb R_+)$.
The model~\eqref{PDV} has a pathwise unique continuous
$(\mathcal F_t)$-adapted solution on $\mathbb R_+$.
Moreover, $\mathbb P$-almost surely, $\bar R_2$ is locally absolutely continuous and satisfies, for almost every $t\ge0$,
\begin{equation}\label{eq:variance-dynamics}
\dot{\bar R}_{2,t}=\dot g_2(t)+K_2(t,t)\sigma_t^2
+\int_0^t\partial_tK_2(s,t)\sigma_s^2\,ds,
\end{equation}
and

\[
\forall t\ge0,\qquad \bar R_{2,t}\ge g_2(t)>0.
\]

\end{prop2}

\begin{proof}
Differentiation under the integral in \eqref{eq:initial-spectral-representation},
justified by~\eqref{eq:initial-spectral-regularity}, gives $g_1\in C^1(\mathbb R_+)$.
The kernels satisfy the estimates used in the proof of
\cite[Theorem~1]{AndresJourdain2023} for $0\le s\le t$.
In that article, the initial inputs are defined differently, by integrals over
the prescribed history, including a Brownian integral over negative times for
$g_1$. However, the same existence and uniqueness argument applies to
deterministic continuous inputs bounded on compact intervals, with $g_2$
bounded away from zero on each compact interval. For every $T>0$, continuity and strict positivity of $g_2$ in Assumption~\ref{assumption3} give $\inf_{t\in[0,T]}g_2(t)>0$.
Fubini's theorem gives~\eqref{eq:variance-dynamics}. In~\eqref{PDV}, $K_2\ge0$ and $g_2>0$ by Assumption~\ref{assumption3}, hence $\bar R_{2,t}\ge g_2(t)>0$.
\end{proof}

\begin{rk}
The integral defining $\bar R_2$ may exist even when $K_2$ is singular on
the diagonal. Our proof uses $K_2(t,t)<\infty$ because
\eqref{eq:variance-dynamics} contains the term $K_2(t,t)\sigma_t^2$.
\end{rk}

\section{Positivity}
\subsection{Positivity and martingality}
The asset-price bound holds under weaker hypotheses than Assumptions~\ref{assumption4} and~\ref{assumption3}: complete monotonicity of $K_1$ and differentiability of the kernels and initial inputs are not required. We first derive the integrated resolvent identity used in its proof.

\begin{lemma}\label{lem:volatility-resolvent}
Suppose that $K_1\in L^2_{\mathrm{loc}}(\mathbb R_+)$ has a resolvent
of the first kind $\mathcal R$.
For $T>0$, suppose that \eqref{PDV} has a
continuous adapted solution on $[0,T]$, with $g_1\in C([0,T])$
and $\bar R_2\ge0$.
Almost surely, the following identity holds:
\begin{equation}\label{eq:resolvent-volatility-integrated}
\forall t\in[0,T],\qquad
\left(\mathcal R*(\sigma-\beta_0-\beta_1g_1-\beta_2\sqrt{\bar R_2})\right)(t)
=\beta_1\int_0^t\sigma_s\,dW_s.
\end{equation}
\end{lemma}

\begin{proof}
Identity~(2.14) in \cite[Lemma~2.6]{AbiJaberLarssonPulido2019}
applies with
\[
X_t\mathrel{:=}\bar R_{1,t}-g_1(t),\qquad X_0=0,\qquad
\mathcal M_t\mathrel{:=}\int_0^t\sigma_s\,dW_s,\qquad K=K_1,\qquad L=\mathcal R.
\]
Indeed, \eqref{PDV} gives $X=K_1*d\mathcal M$, and both $X$ and $\sigma$
are continuous and adapted. Hence, almost surely, we have
\[
\forall t\in[0,T],\qquad
\left(\mathcal R*(\bar R_1-g_1)\right)(t)
=\int_0^t\sigma_s\,dW_s.
\]
Multiplying by $\beta_1$ and using
$\sigma-\beta_0-\beta_1g_1-\beta_2\sqrt{\bar R_2}
=\beta_1(\bar R_1-g_1)$ proves
\eqref{eq:resolvent-volatility-integrated}.
\end{proof}

\begin{thm2}\label{thm:GL_positive_moments_compact}
Let $T>0$. Assume that $g_1,g_2\in C([0,T])$, with $g_2\ge0$, and that
\eqref{PDV} has a unique continuous $(\mathcal F_t)$-adapted solution
on $[0,T]$ under $\mathbb P$.
Assume that $\sigma_t\ge0$ for every $t\in[0,T]$,
$\mathbb P$-almost surely.
Suppose that $K_2$ is nonnegative and measurable and that
$K_1$ is a nonnegative convolution
kernel in $L^2_{\mathrm{loc}}(\mathbb R_+)$, with a resolvent
of the first kind $\mathcal R$, and that
\begin{equation}\label{eq:C2T}
C_{2,T}:=
\sup_{r\in[0,T]}\sup_{u\in[0,r]}
\int_{[0,r-u]}K_2(u,r-s)\,\mathcal R(ds)<\infty.
\end{equation}
Then, $\mathbb P$-almost surely, the following bound holds:
\begin{align}
\label{eq:GL_asset_pathwise_upper_bound}
\sup_{t\in[0,T]}S_t
&\le S_0\exp\left(
\mathcal R([0,T])
\left(
\frac{\beta_0}{-\beta_1}
+\sup_{t\in[0,T]}(-g_1(t))_+
+\frac{\beta_2}{-\beta_1}
\sqrt{\sup_{t\in[0,T]}g_2(t)}
\right)
+\frac{\beta_2^2}{2\beta_1^2}
\mathcal R([0,T])C_{2,T}
\right).
\end{align}
Hence, under $\mathbb P$, $S$ is a bounded martingale
on $[0,T]$.
\end{thm2}

\begin{rk}\label{rk:asset-bound-intuition}
Under Assumptions~\ref{assumption4} and~\ref{assumption3}, the identity~\eqref{eq:resolvent-volatility-integrated} explains why nonnegative volatility limits how large $\int_0^t\sigma_s\,dW_s$ can be. Indeed, since $\beta_1<0$, an increase in this integral pushes $\sigma$ downward. The Brownian integral is therefore controlled and, since it drives the logarithm of the asset price, this yields an upper bound on the price.
\end{rk}

\begin{rk}\label{rk:vanishing-right-wing}
Under the assumptions of Theorem~\ref{thm:GL_positive_moments_compact},
for fixed deterministic initial inputs, let $B_T$ denote the
right-hand side of \eqref{eq:GL_asset_pathwise_upper_bound}.
Then, for every $p>0$,
\[
\E\!\left[\sup_{t\in[0,T]}S_t^p\right]\le B_T^p<\infty.
\]
European call options with maturity $T$ and strike $k\ge B_T$
have zero payoff almost surely, and hence zero price:
$\E[(S_T-k)_+]=0$.
The Black--Scholes implied volatility smile is therefore
identically zero for strikes $k\ge B_T$.
Theorem~\ref{thm:GL_positive_moments_compact} provides a theoretical basis for the vanishing right wing of the implied volatility smile observed numerically by Guyon and Lekeufack \cite{GuyonLekeufack2023} and by Gazzani and Guyon \cite{GazzaniGuyon2025}. This vanishing behaviour motivated the quadratic correction introduced by Guyon and Lekeufack \cite[Section~4.5]{GuyonLekeufack2023} to reproduce the rising right wing of the SPX smile. Alternatively, Guyon and Parent \cite{guyon2025discrete} have shown that a discrete-time version of the model with fat-tailed innovations captures market smiles very well, even for large strikes.

\end{rk}

\begin{rk}
If $K_1$ is continuous at zero and $K_1(0)>0$, then
$\mathcal R(\{0\})=\frac{1}{K_1(0)}$ by $K_1*\mathcal R=1$. Definition~\eqref{eq:C2T} gives
\[
\frac{1}{K_1(0)}
\sup_{r\in[0,T]}\sup_{u\in[0,r]}K_2(u,r)
\le C_{2,T}
\le
\mathcal R([0,T])
\sup_{r\in[0,T]}\sup_{u\in[0,r]}K_2(u,r).
\]
Then condition~\eqref{eq:C2T} holds when $K_2$ is bounded on
$\{(u,r):0\le u\le r\le T\}$, for example when
$K_2(s,t)=K_2(0,t-s)$ and
$u\mapsto K_2(0,u)$ belongs to $C(\mathbb R_+)$.
\end{rk}

\begin{rk}
The H\"older regularity conditions
\begin{align}
\left(
\int_0^{T-h}|K_1(u+h)-K_1(u)|^2\,du
+\int_0^hK_1(u)^2\,du
\right)^{\frac{1}{2}}
&\le C_Th^{\delta_1},\label{eq:K1-holder-increment}\\
\sup_{t\in[0,T-h]}\left(
\int_0^t|K_2(s,t+h)-K_2(s,t)|\,ds
+\int_t^{t+h}K_2(s,t+h)\,ds
\right)
&\le C_Th^{\delta_2},\label{eq:K2-holder-increment}
\end{align}
for every $T>0$ and $h\in[0,T]$, with a constant $C_T>0$,
follow from $K_1\in C^1(\mathbb R_+)$ in Notation~\ref{not:return-kernel} and~\eqref{eq:K2-local-bounds}, with $\delta_1=\frac{1}{2}$ and $\delta_2=1$.
Under the increment conditions~\eqref{eq:K1-holder-increment} and~\eqref{eq:K2-holder-increment}, \cite[Theorem~1~(B.i)]{GuyonJeanninJourdain2026Martingale} shows that nonnegative volatility implies martingality.
(For simplicity, the result is stated there for a convolution kernel $K_2$,
but the proofs remain unchanged.)
Theorem~\ref{thm:GL_positive_moments_compact} of the present paper gives moreover a deterministic pathwise bound on $S$ over $[0,T]$, under an additional condition on the return kernel $K_1$: it needs to admit a resolvent of the first kind. Assumption~\ref{assumption4} ensures this, while Assumption~\ref{assumption3} implies~\eqref{eq:C2T}.
\end{rk}

\begin{proof}[Proof of Theorem~\ref{thm:GL_positive_moments_compact}]
By Lemma~\ref{lem:volatility-resolvent}
and $\sigma-\beta_0-\beta_1g_1-\beta_2\sqrt{\bar R_2}=\beta_1(\bar R_1-g_1)$,
almost surely, we have
\[
\forall t\in[0,T],\qquad
\int_{[0,t]}(\bar R_{1,t-s}-g_1(t-s))\,\mathcal R(ds)
=\int_0^t\sigma_s\,dW_s.
\]
Moreover, $\sigma\ge0$ and $\beta_1<0$ imply
\begin{align*}
\forall v\in[0,T],\quad
\bar R_{1,v}-g_1(v)
&=\frac{\beta_0+\beta_2\sqrt{\bar R_{2,v}}-\sigma_v}{-\beta_1}
-g_1(v) \le\frac{\beta_0}{-\beta_1}
+\frac{\beta_2}{-\beta_1}\sqrt{\bar R_{2,v}}
+\sup_{u\in[0,T]}(-g_1(u))_+.
\end{align*}
Since $S_0>0$, the stochastic exponential solution of~\eqref{PDV} satisfies $S_t>0$, so $\log\left(\frac{S_t}{S_0}\right)$ is well defined. It\^o's formula, the nonnegativity of $\mathcal R$, and
$\mathcal R([0,t])\le\mathcal R([0,T])$ therefore give
\begin{align}
\log\frac{S_t}{S_0}
&=\int_0^t\sigma_u\,dW_u-\frac{1}{2}\int_0^t\sigma_u^2\,du
=\int_{[0,t]}(\bar R_{1,t-s}-g_1(t-s))\,\mathcal R(ds)
-\frac{1}{2}\int_0^t\sigma_u^2\,du
\label{eq:direct-resolvent-identity}\\
&\le\mathcal R([0,T])
\left(\frac{\beta_0}{-\beta_1}
+\sup_{u\in[0,T]}(-g_1(u))_+\right)
+\frac{\beta_2}{-\beta_1}
\int_{[0,t]}\sqrt{\bar R_{2,t-s}}\,\mathcal R(ds)-\frac{1}{2}\int_0^t\sigma_u^2\,du.\label{eq:asset-resolvent-upper-bound}
\end{align}
By definition of ~\eqref{PDV},
\[
\bar R_{2,t}=g_2(t)+\int_0^t K_2(s,t)\sigma_s^2\,ds.
\]
The Cauchy--Schwarz inequality, Tonelli's theorem, and
the definition~\eqref{eq:C2T} of $C_{2,T}$ give
\begin{align}
\left(\int_{[0,t]}\sqrt{\bar R_{2,t-s}}\,\mathcal R(ds)\right)^2
&\le\mathcal R([0,T])\int_{[0,t]}\bar R_{2,t-s}\,\mathcal R(ds)
\le\mathcal R([0,T])^2\sup_{u\in[0,T]}g_2(u)\notag\\
&\quad+\mathcal R([0,T])\int_{u\in[0,t]}\sigma_u^2
\int_{s\in[0,t-u]}K_2(u,t-s)\,\mathcal R(ds)\,du\notag\\
&\le\mathcal R([0,T])^2\sup_{u\in[0,T]}g_2(u)+\mathcal R([0,T])C_{2,T}\int_0^t\sigma_u^2\,du.\label{eq:asset-resolvent-cauchy-bound}
\end{align}
Combining~\eqref{eq:asset-resolvent-cauchy-bound} with~\eqref{eq:asset-resolvent-upper-bound}
and the inequalities $\sqrt{x+y}\le\sqrt x+\sqrt y$ and
$ax-\frac{x^2}{2}\le \frac{a^2}{2}$ gives
\begin{align*}
\forall t\in[0,T],\quad
\log\frac{S_t}{S_0}
&\le\mathcal R([0,T])
\left(\frac{\beta_0}{-\beta_1}
+\sup_{u\in[0,T]}(-g_1(u))_+
+\frac{\beta_2}{-\beta_1}\sqrt{\sup_{u\in[0,T]}g_2(u)}\right)\\
&\quad+\frac{\beta_2}{-\beta_1}
\sqrt{\mathcal R([0,T])C_{2,T}}\sqrt{\int_0^t\sigma_u^2\,du}
-\frac{1}{2}\int_0^t\sigma_u^2\,du\\
&\le\mathcal R([0,T])
\left(\frac{\beta_0}{-\beta_1}
+\sup_{u\in[0,T]}(-g_1(u))_+
+\frac{\beta_2}{-\beta_1}\sqrt{\sup_{u\in[0,T]}g_2(u)}\right)
\\
&\quad+\frac{\beta_2^2}{2\beta_1^2}\mathcal R([0,T])C_{2,T}.
\end{align*}
Exponentiation gives~\eqref{eq:GL_asset_pathwise_upper_bound}.
The upper bound is deterministic.
Thus $S$ is a bounded local martingale, hence a martingale.
\end{proof}

\noindent \textbf{Extensions to stochastic volatility models.}\label{par:price-bound-extensions}
The following proposition is a direct adaptation of the proof of Theorem~\ref{thm:GL_positive_moments_compact}. It applies to several stochastic volatility models at perfect negative correlation between their Brownian drivers, $\rho=-1$, as illustrated in Table~\ref{tab:general-price-bound-models}.

\begin{prop2}\label{prop:general-volterra-price-bound}
Let $T>0$, $S_0,c>0$, and let $\sigma$ and $G$ be continuous adapted processes on $[0,T]$ satisfying
\begin{equation}\label{eq:general-volterra-volatility}
\forall t\in[0,T],\qquad
\begin{aligned}
\frac{dS_t}{S_t}&=\sigma_t\,dW_t,\\
\Phi(\sigma_t)&=G_t-c\int_0^tK(t-s)\sigma_s\,dW_s,
\end{aligned}
\end{equation}
where $\Phi:\mathbb R\to\mathbb R$ is continuous and $\Phi(\sigma_t)\ge0$ almost surely for every $t\in[0,T]$.
Suppose that $K\in L^2_{\mathrm{loc}}(\mathbb R_+)$ admits a nonnegative resolvent of the first kind $\mathcal R$. If deterministic constants $a_T\in\mathbb R$, $b_T\ge0$ and $d_T<\frac{c}{2}$ satisfy, almost surely,
\begin{equation}\label{eq:general-volterra-input-bound}
\forall t\in[0,T],\qquad
(\mathcal R*G)(t)\le a_T+b_T\sqrt{\int_0^t\sigma_s^2\,ds}
+d_T\int_0^t\sigma_s^2\,ds,
\end{equation}
then, almost surely,
\begin{equation}\label{eq:general-volterra-price-bound}
\sup_{t\in[0,T]}S_t\le S_0\exp\left(
\frac{a_T}{c}+\frac{b_T^2}{2c(c-2d_T)}\right).
\end{equation}
In particular, $S$ is a bounded martingale on $[0,T]$.
\end{prop2}

\begin{proof}
Convolution with $\mathcal R$, justified as in Lemma~\ref{lem:volatility-resolvent}, gives
\begin{align*}
\log\frac{S_t}{S_0}
&=\frac{(\mathcal R*G)(t)-(\mathcal R*\Phi(\sigma))(t)}{c}
-\frac{1}{2}\int_0^t\sigma_s^2\,ds\\
&\le\frac{a_T}{c}+\frac{b_T}{c}\sqrt{\int_0^t\sigma_s^2\,ds}
-\frac{c-2d_T}{2c}\int_0^t\sigma_s^2\,ds\\
&\le\frac{a_T}{c}+\frac{b_T^2}{2c(c-2d_T)}.
\end{align*}
The deterministic bound makes the local martingale $S$ a martingale.
\end{proof}
\begin{example}\label{ex:general-price-bound-models}
Under the hypotheses of Theorem~\ref{thm:GL_positive_moments_compact}, \eqref{PDV} has the form~\eqref{eq:general-volterra-volatility} with $\Phi(x)=x$, $c=-\beta_1$, $K=K_1$ and $G=\beta_0+\beta_1g_1+\beta_2\sqrt{\bar R_2}$. Estimate~\eqref{eq:asset-resolvent-cauchy-bound} gives~\eqref{eq:general-volterra-input-bound} with $d_T=0$, recovering~\eqref{eq:GL_asset_pathwise_upper_bound}.

For the Heston and rough Heston models, $\Phi(x)=x^2$, $c=\nu$ and
\[
G_t=\sigma_0^2+\int_0^tK(t-s)\kappa(\theta-\sigma_s^2)\,ds,
\qquad
(a_T,b_T,d_T)=\left(\sigma_0^2\mathcal R([0,T])+\kappa\theta T,0,-\kappa\right).
\]
For the Heston model, we have $K=1$ and $\mathcal R=\delta_0$, recovering the bound in \cite[Proposition~4.2.8]{Jacquier2017Notes}. For the rough Heston model, we have
\[
K(t)=\frac{t^{\alpha-1}}{\Gamma(\alpha)},\quad
\alpha\in\left(\frac{1}{2},1\right),\qquad
\mathcal R([0,T])=\frac{T^{1-\alpha}}{\Gamma(2-\alpha)}.
\]

Table~\ref{tab:general-price-bound-models} collects the resulting price bounds at $\rho=-1$.

\begin{table}[H]
\centering
\small
\renewcommand{\arraystretch}{1.5}
\begin{tabularx}{\textwidth}{@{}X p{0.35\textwidth}@{}}
\toprule
Model and volatility dynamics at $\rho=-1$ & Upper bound for
\newline $\log\left(\frac{S_t}{S_0}\right)$, $0\le t\le T$ \\
\midrule
\textit{Heston model} \cite{Heston1993}
\newline $d(\sigma_t^2)=\kappa(\theta-\sigma_t^2)\,dt-\nu\sigma_t\,dW_t$
& $\displaystyle \frac{\sigma_0^2+\kappa\theta T}{\nu}$ \\
\midrule
\textit{Rough Heston model} \cite{ElEuchRosenbaum2019}
\newline $\displaystyle\begin{aligned}[t]
\sigma_t^2&=\sigma_0^2+\int_0^tK(t-s)\kappa(\theta-\sigma_s^2)\,ds\\
&\quad-\nu\int_0^tK(t-s)\sigma_s\,dW_s
\end{aligned}$
& $\displaystyle \frac{\sigma_0^2T^{1-\alpha}}{\nu\Gamma(2-\alpha)}+\frac{\kappa\theta T}{\nu}$ \\
\midrule
\textit{Lognormal SABR model} ($\beta=1$) \cite{HaganKumarLesniewskiWoodward2002}
\newline $d\sigma_t=-\nu\sigma_t\,dW_t$
& $\displaystyle \frac{\sigma_0}{\nu}$ \\
\midrule
\textit{Sch\"obel--Zhu model} \cite{SchobelZhu1999}
\newline $d\sigma_t=\kappa(\theta-\sigma_t)\,dt-\nu\,dW_t$
& $\displaystyle \frac{\sigma_0^2}{2\nu}+\frac{\nu T}{2}+\frac{\kappa^2\theta^2T}{2\nu(2\kappa+\nu)}$ \\
\bottomrule
\end{tabularx}
\caption{Applications of Proposition~\ref{prop:general-volterra-price-bound} with deterministic initial values and $\kappa,\theta,\nu>0$.}
\label{tab:general-price-bound-models}
\end{table}
Thus, the endogenous specifications in Table~\ref{tab:general-price-bound-models} at $\rho=-1$ have an implied volatility smile with a flat, zero right wing beyond a finite strike. Reproducing a rising right wing motivates an additional volatility driver independent of $W$ or, in the PDV model, a volatility function that increases again for sufficiently large positive values of the past-return factor, as discussed in Remark~\ref{rk:vanishing-right-wing}.
\end{example}

\subsection{Sufficient conditions for positivity}
We now return to the volatility $\sigma$ in \eqref{PDV}.

\begin{notation}\label{not:initial-cumulative}
Under Assumption~\ref{assumption4}, we set
\[
\underline\lambda_1\mathrel{:=}\inf\operatorname{supp}\nu.
\]
We define the nonnegative measure
\[
\hat\nu(d\lambda)\mathrel{:=}\frac{\lambda(\lambda-\underline\lambda_1)}{\Lambda_1}
\,\nu(d\lambda).
\]
By condition~\eqref{eq:nu-second-moment}, it is finite and has total mass $\lambda_{1,\mathrm{eff}}-\underline\lambda_1$.
For $\lambda>0$, we define
\[
\eta_0(\lambda)\mathrel{:=}\beta_1\frac{\Lambda_1}{\lambda}R_{1,0}(\lambda)-\sigma_0,
\qquad
D\mathrel{:=}\sup_{r\in[\underline\lambda_1,\infty)}
\int_{\lambda\in(0,r]}\eta_0(\lambda)\,\hat\nu(d\lambda).
\]
Conditions~\eqref{eq:nu-second-moment} and~\eqref{eq:initial-spectral-regularity} give
\begin{equation}\label{eq:initial-d-integrability}
\int_{\lambda\in(0,\infty)}|\eta_0(\lambda)|\,\hat\nu(d\lambda)
\le|\beta_1|\int_{\lambda\in(0,\infty)}\lambda|R_{1,0}(\lambda)|\,\nu(d\lambda)
+|\sigma_0|(\lambda_{1,\mathrm{eff}}-\underline\lambda_1)<\infty.
\end{equation}
Hence $\eta_0\in L^1(\hat\nu)$ and
$D\le\int_{\lambda\in(0,\infty)}|\eta_0(\lambda)|\,\hat\nu(d\lambda)<\infty$.
Moreover, $\hat\nu(\{\underline\lambda_1\})=0$ gives $D\ge0$.
\end{notation}

Under Assumption~\ref{assumption4}, $\mathcal R-\Lambda_1^{-1}\delta_0$ has a locally integrable density. Differentiating $K_1*\mathcal R=1$ with
\cite[Chapter~3, Corollary~7.4(ii)]{GripenbergLondenStaffans1990} gives
\begin{equation}\label{eq:resolvent-measure-identity}
\Lambda_1\mathcal R(dt)=\delta_0(dt)-(K_1'*\mathcal R)(t)\,dt.
\end{equation}
Since $\mathcal R(\{0\})=\Lambda_1^{-1}$, we have
\[
(K_1'*\mathcal R)(t)=\frac{K_1'(t)}{\Lambda_1}
+\int_{(0,t]}K_1'(t-s)\,\mathcal R(ds).
\]
Continuity of $K_1'$ and local finiteness of $\mathcal R$ imply that
the integral tends to zero as $t\downarrow0$. Consequently,
\begin{equation}\label{eq:resolvent-memory-initial}
(K_1'*\mathcal R)(0)=\frac{K_1'(0)}{\Lambda_1}
=-\lambda_{1,\mathrm{eff}}.
\end{equation}
By \cite[Chapter~5, Theorem~5.4]{GripenbergLondenStaffans1990},
the density of $\mathcal R-\Lambda_1^{-1}\delta_0$ is completely monotone.
Identity~\eqref{eq:resolvent-measure-identity} therefore shows that
$-(K_1'*\mathcal R)$ is completely monotone on $(0,\infty)$.
Complete monotonicity and~\eqref{eq:resolvent-memory-initial} imply local absolute continuity on $\mathbb R_+$.
Thus $K_1'*\mathcal R$ is nonpositive and nondecreasing, which gives
\begin{equation}\label{eq:resolvent-memory-bound}
\begin{aligned}
&\text{for almost every }t>0,\quad (K_1'*\mathcal R)'(t)\ge0,\\
&\forall t\ge0,\quad
\int_0^t(K_1'*\mathcal R)'(s)\,ds
=(K_1'*\mathcal R)(t)+\lambda_{1,\mathrm{eff}}
\le\lambda_{1,\mathrm{eff}}.
\end{aligned}
\end{equation}

\begin{lemma}\label{lem:volatility-semimartingale}
Under Assumptions~\ref{assumption4} and~\ref{assumption3}, for every $T>0$,
the process $\sigma$ is a semimartingale with dynamics
\begin{align}
d\sigma_t
&=\left(-\lambda_{1,\mathrm{eff}}\sigma_t
+((K_1'*\mathcal R)'*\sigma)(t)\right.\notag\\*
&\qquad\left.+\Lambda_1\frac{d}{dt}
\left(\mathcal R*(\beta_0+\beta_1g_1+\beta_2\sqrt{\bar R_2})\right)(t)
\right)\,dt+\beta_1\Lambda_1\sigma_t\,dW_t,
\label{eq:resolvent-sigma-dynamics}
\end{align}
where the drift is defined for almost every $t\in[0,T]$, almost surely.
\end{lemma}

\begin{proof}
By Proposition~\ref{prop:existence-uniqueness-moments},
almost surely, $g_1$ and $\bar R_2$ are absolutely continuous on $[0,T]$
and $\inf_{t\in[0,T]}\bar R_{2,t}\ge\min_{t\in[0,T]}g_2(t)>0$ by Assumption~\ref{assumption3}.
The chain rule therefore gives absolute continuity of $\sqrt{\bar R_2}$.
Convolving \eqref{eq:resolvent-measure-identity} with $\sigma$ and using
\eqref{eq:resolvent-volatility-integrated} gives
\begin{align}
\forall t\in[0,T],\qquad
\sigma_t
&=\Lambda_1
\left(\mathcal R*(\beta_0+\beta_1g_1+\beta_2\sqrt{\bar R_2})\right)(t)
+((K_1'*\mathcal R)*\sigma)(t)
+\beta_1\Lambda_1\int_0^t\sigma_s\,dW_s.
\label{eq:resolvent-sigma-integrated}
\end{align}
We show that the two convolution terms in \eqref{eq:resolvent-sigma-integrated}
are absolutely continuous, so that their time derivatives can be taken
to obtain \eqref{eq:resolvent-sigma-dynamics}.
By~\eqref{eq:resolvent-memory-initial} and~\eqref{eq:resolvent-memory-bound},
$K_1'*\mathcal R$ is absolutely continuous on $[0,T]$ and
$(K_1'*\mathcal R)(0)=-\lambda_{1,\mathrm{eff}}$. Hence
\begin{align}\label{eq:resolvent-memory-convolution-derivative}
\frac{d}{dt}\left((K_1'*\mathcal R)*\sigma\right)(t)
=(K_1'*\mathcal R)(0)\sigma_t+((K_1'*\mathcal R)'*\sigma)(t)
=-\lambda_{1,\mathrm{eff}}\sigma_t+((K_1'*\mathcal R)'*\sigma)(t).
\end{align}
The right-hand side is integrable on $[0,T]$, so
$(K_1'*\mathcal R)*\sigma$ is absolutely continuous.
For the other convolution in \eqref{eq:resolvent-sigma-integrated},
the measure identity~\eqref{eq:resolvent-measure-identity} gives
\begin{align}\label{eq:resolvent-input-convolution-identity}
\Lambda_1\left(\mathcal R*(\beta_0+\beta_1g_1+\beta_2\sqrt{\bar R_2})\right)
=\beta_0+\beta_1g_1+\beta_2\sqrt{\bar R_2}-(K_1'*\mathcal R)*(\beta_0+\beta_1g_1+\beta_2\sqrt{\bar R_2}).
\end{align}
The function $\beta_0+\beta_1g_1+\beta_2\sqrt{\bar R_2}$ is absolutely continuous.
Formula~\eqref{eq:resolvent-memory-convolution-derivative}, applied with
$\sigma$ replaced by this continuous function, shows that its convolution
with $K_1'*\mathcal R$ is absolutely continuous as well.
Identity~\eqref{eq:resolvent-input-convolution-identity} therefore shows that the first convolution in
\eqref{eq:resolvent-sigma-integrated} is absolutely continuous, completing the proof.
\end{proof}

\begin{assumption}\label{assumption:positivity}
The following conditions hold:
\begin{enumerate}[(i)]
\myitem[(i)]{item:positivity-initial-volatility} $\sigma_0\ge0$.
\myitem[(ii)]{item:positivity-variance} For $a\in\mathbb R$, we consider, for each $s\ge0$ and almost every $t\ge0$,
\begin{align}
\text{for }t\ge s,\qquad
\partial_tK_2(s,t)+(2\underline\lambda_1-a)K_2(s,t)
&\ge0,\label{not:def-m}\\
\dot g_2(t)+(2\underline\lambda_1-a)g_2(t)
&\ge0.\label{eq:initial-variance-growth}
\end{align}
We define the set
\[
\left\{a\in\mathbb R:
\eqref{not:def-m}\text{ and }\eqref{eq:initial-variance-growth}\text{ hold}\right\}.
\]
We assume that this set is nonempty. It is closed and, by integration of~\eqref{eq:initial-variance-growth} over $[0,1]$, bounded above by $2\underline\lambda_1+\log\frac{g_2(1)}{g_2(0)}$, where $g_2(0),g_2(1)>0$ by Assumption~\ref{assumption3}. We define $m$ as its maximum.
We also assume that either $m\ge0$, or $m<0$ and
\[
\forall t\ge0,\qquad
\begin{aligned}
&g_2(t)\le g_2(0),\\
&\text{for almost every }s\ge t,\quad
\partial_sK_2(t,s)+2\underline\lambda_1K_2(t,s)\ge0.
\end{aligned}
\]
\myitem[(iii)]{item:positivity-initial-balance} The following bound holds:
\begin{equation}\label{Dm}
D\le\beta_0\underline\lambda_1+\frac{\beta_2m}{2}\sqrt{\bar R_{2,0}}.
\end{equation}
\end{enumerate}
\end{assumption}

\begin{rk}\label{rk:positivity-assumption-roles}
By~\eqref{eq:resolvent-sigma-dynamics}, the diffusion coefficient $\beta_1\Lambda_1\sigma_t$ vanishes at $\sigma_t=0$, as does the drift term $-\lambda_{1,\mathrm{eff}}\sigma_t$. By the nonnegativity of $(K_1'*\mathcal R)'$ in~\eqref{eq:resolvent-memory-bound}, $((K_1'*\mathcal R)'*\sigma)(t)\ge0$ if $\sigma_s\ge0$ for all $s\in[0,t]$. To prevent the drift from pushing volatility below zero at $\sigma_t=0$, it therefore remains to show that
\[
\text{for almost every }t\ge0,\qquad
\Lambda_1\frac{d}{dt}
\left(\mathcal R*\left(\beta_0+\beta_1g_1+\beta_2\sqrt{\bar R_2}\right)\right)(t)
\ge0.
\]
The assumptions control the contributions of $\sqrt{\bar R_2}$ and $g_1$ to this term as follows.
\begin{itemize}
\item For a single exponential in each kernel, with initial inputs as in~\eqref{eq:exponential-inputs}, at $\sigma=0$, $\beta_1\bar R_1$ relaxes toward zero at rate $\lambda_{1,1}$, raising $\sigma$, whereas $\beta_2\sqrt{\bar R_2}$ decays at rate $\frac{\lambda_{2,1}}{2}$, lowering $\sigma$. Similarly, Condition~\myref{item:positivity-variance} controls the contribution of $\sqrt{\bar R_2}$ by limiting its decay relative to the return kernel $K_1$. More precisely,  \eqref{eq:initial-variance-growth} controls the decay of $g_2(t)$, while~\eqref{not:def-m} controls that of the integral term $\int_0^tK_2(s,t)\sigma_s^2\,ds$.

\item Condition~\myref{item:positivity-initial-balance} bounds, through $D$, the potentially negative drift contribution of $\beta_1g_1$.
\end{itemize}
\end{rk}

\begin{rk}\label{rk:variance-logarithmic-decay}
Since $g_2>0$ by Assumption~\ref{assumption3}, the definition of $m$ in Assumption~\ref{assumption:positivity}~\myref{item:positivity-variance} equivalently reads
\[
m=2\underline\lambda_1-\max\left\{
\operatorname*{ess\,sup}_{t\ge0}\frac{-\dot g_2(t)}{g_2(t)},\quad
\sup_{s\ge0}\operatorname*{ess\,sup}_{\substack{t\ge s\\K_2(s,t)>0}}
\frac{-\partial_tK_2(s,t)}{K_2(s,t)}
\right\}.
\]
\end{rk}

We define $Z$ by
\begin{equation}\label{eq:def-Z}
dZ_t=-\lambda_{1,\mathrm{eff}}Z_t\,dt+\beta_1\Lambda_1Z_t\,dW_t,
\qquad Z_0=1,
\end{equation}
i.e., integrating this SDE,
\[
\forall t\ge0,\qquad
Z_t=\exp\left(\beta_1\Lambda_1W_t-\left(\lambda_{1,\mathrm{eff}}+\frac{\beta_1^2\Lambda_1^2}{2}\right)t\right)>0.
\]

\begin{thm2}\label{prop:pos-MR1-generalK2-vector0}
Under Assumptions~\ref{assumption4}, \ref{assumption3}, and~\ref{assumption:positivity}, $\mathbb P$-almost surely, the following bound holds:
\begin{equation}\label{eq:positivity-pathwise-bound}
\forall t\ge0,\qquad
\sigma_t\ge Z_t\left(\sigma_0+
\left(\beta_0\underline\lambda_1
-\left(D-\frac{\beta_2m}{2}\sqrt{\bar R_{2,0}}\right)_+\right)
\int_0^tZ_s^{-1}\,ds\right)\ge0.
\end{equation}
\begin{itemize}
\item If $\sigma_0>0$, the lower bound is strictly positive for every $t\ge0$.

\item If $\sigma_0=0$, it is strictly positive for every $t>0$ provided that
$$
\beta_0\underline\lambda_1>
\left(D-\frac{\beta_2m}{2}\sqrt{\bar R_{2,0}}\right)_+.
$$
\end{itemize}
Moreover, for every $T>0$, $(S_t)_{t\in[0,T]}$ is a bounded martingale
under $\mathbb P$.
\end{thm2}

\begin{rk}
The positivity proof requires $K_1(0)<\infty$ and therefore
does not cover kernels with $K_1(t)\to\infty$ as $t\downarrow0$.
Whether volatility reaches zero can depend on whether the kernel is
regular or rough.
For Volterra CIR processes with regular kernels,
Bondi and Pulido~\cite[Theorem~9~(ii)]{BondiPulido2024} give a Feller-type condition
ensuring that zero is inaccessible.
By contrast, the variance in the fractional rough Heston model has a positive
atom at zero at every positive time, see
\cite[Theorem~3]{BourdonJeannin2026}
and \cite[Theorem~1.6, equation~(1.9)]{FriesenGerholdWiedermann2026}.
\end{rk}

\begin{rk}
For a diffusion
\[
dY_t=b(Y_t)\,dt+\Sigma(Y_t)\,dW_t,
\qquad \mathcal C\mathrel{:=}\{y\in\mathbb R^d:y_1\ge0\},
\]
\cite[Theorem~4.1]{FilipovicTappeTeichmann2014} characterizes
local stochastic invariance of $\mathcal C$, meaning that solutions starting
in $\mathcal C$ remain there locally. Under the hypotheses of
\cite[Theorem~4.1]{FilipovicTappeTeichmann2014}, the criterion is given by
\[
\forall y\in\mathbb R^d\text{ with }y_1=0,\qquad
\Sigma_{1,\cdot}(y)=0,\qquad b_1(y)\ge0.
\]
In~\eqref{eq:resolvent-sigma-dynamics}, the coefficient of $dW_t$ vanishes at $\sigma_t=0$. The drift depends on the past through $((K_1'*\mathcal R)'*\sigma)(t)$ and on $\bar R_{2,t}$, so $\sigma$ alone is not a diffusion with coefficients determined by its current value. The proof applies Gr\"onwall's lemma to $\sigma^-$ without invoking this diffusion criterion.
\end{rk}

\begin{proof}[Proof of Theorem~\ref{prop:pos-MR1-generalK2-vector0}]
Lemma~\ref{lem:volatility-semimartingale}
gives the decomposition~\eqref{eq:resolvent-sigma-dynamics}. The dynamics~\eqref{eq:variance-dynamics} and conditions~\eqref{not:def-m} and~\eqref{eq:initial-variance-growth} give
$\dot{\bar R}_{2,t}+(2\underline\lambda_1-m)\bar R_{2,t}\ge0$ for almost every $t\ge0$. Since $\bar R_{2,0}=g_2(0)>0$ by~\eqref{eq:initial-values},
\begin{equation}\label{eq:variance-lower-bound}
\forall t\ge0,\qquad \bar R_{2,t}\ge e^{-(2\underline\lambda_1-m)t}\bar R_{2,0}>0.
\end{equation}
In~\eqref{eq:resolvent-sigma-dynamics}, when $\sigma_t=0$, both $\beta_1\Lambda_1\sigma_t$ and $-\lambda_{1,\mathrm{eff}}\sigma_t$ vanish.

\medskip

\noindent\textit{Step 1: Nonnegativity of the last drift term in \eqref{eq:resolvent-sigma-dynamics}.} Since $\mathcal R*K_1=1$, we have
\begin{align*}
\Lambda_1\frac{d}{dt}
\left(\mathcal R*(\beta_0+\beta_1g_1+\beta_2\sqrt{\bar R_2})\right)(t)
&=\Lambda_1\frac{d}{dt}
\left[\mathcal R*
\left(\beta_0+\beta_1g_1+\beta_2\sqrt{\bar R_2}
-\frac{\sigma_0}{\Lambda_1}K_1\right)\right](t).
\end{align*}
We define, for almost every $t\ge0$,
\begin{equation}\label{eq:positivity-H-definition}
H_t:=\left(\frac{d}{dt}+\underline\lambda_1\right)
\left(\beta_0+\beta_1g_1(t)+\beta_2\sqrt{\bar R_{2,t}}
-\frac{\sigma_0}{\Lambda_1}K_1(t)\right).
\end{equation}
Since $K_1(0)=\Lambda_1$ and
$\sigma_0=\beta_0+\beta_1g_1(0)+\beta_2\sqrt{\bar R_{2,0}}$ by~\eqref{eq:initial-values},
the expression in parentheses in \eqref{eq:positivity-H-definition}
is zero at $t=0$. Hence we have
\[
\beta_0+\beta_1g_1(t)+\beta_2\sqrt{\bar R_{2,t}}
-\frac{\sigma_0}{\Lambda_1}K_1(t)
=\int_0^t e^{-\underline\lambda_1(t-s)}H_s\,ds.
\]
We define, for almost every $t\ge0$,
\[
J(t):=\Lambda_1\frac{d}{dt}
\left(\mathcal R*e^{-\underline\lambda_1\cdot}\right)(t).
\]
Identity~\eqref{eq:resolvent-measure-identity} gives $\Lambda_1\mathcal R(\{0\})=1$.
Convolving with $\mathcal R$ and differentiating gives, for almost every $t\ge0$,
\begin{align}
\Lambda_1\frac{d}{dt}
\left(\mathcal R*(\beta_0+\beta_1g_1+\beta_2\sqrt{\bar R_2})\right)(t)
&=\Lambda_1\frac{d}{dt}
\left[(\mathcal R*e^{-\underline\lambda_1\cdot})*H\right](t)\notag\\
&=\Lambda_1
(\mathcal R*e^{-\underline\lambda_1\cdot})(0)H_t
+\Lambda_1\int_0^t
(\mathcal R*e^{-\underline\lambda_1\cdot})'(t-s)H_s\,ds\notag\\
&=\underbrace{\Lambda_1\mathcal R(\{0\})}_{=\,1}H_t
+\int_0^t J(t-s)H_s\,ds
=H_t+(J*H)(t).
\label{eq:resolvent-H}
\end{align}
We prove that $J\ge0$ and $H\ge0$, which makes the right-hand side of~\eqref{eq:resolvent-H} nonnegative.
Since $K_1*\mathcal R=1$, associativity gives
\begin{align*}
\left(K_1*(\mathcal R*e^{-\underline\lambda_1\cdot})\right)(t)
&=\left((K_1*\mathcal R)*e^{-\underline\lambda_1\cdot}\right)(t)
=\int_0^t e^{-\underline\lambda_1s}\,ds.
\end{align*}
Since $(\mathcal R*e^{-\underline\lambda_1\cdot})(0)=\frac{1}{\Lambda_1}$
and $(\mathcal R*e^{-\underline\lambda_1\cdot})'=\frac{J}{\Lambda_1}$,
differentiation yields
\begin{equation}\label{eq:resolvent-weights}
K_1+J*K_1=\Lambda_1e^{-\underline\lambda_1\cdot}.
\end{equation}
The signed measure $\delta_0+e^{\underline\lambda_1t}J(t)\,dt$ satisfies the resolvent identity for
\[
\frac{e^{\underline\lambda_1t}K_1(t)}{\Lambda_1}
=\frac{1}{\Lambda_1}\int_{\lambda\in(0,\infty)}
\lambda e^{-(\lambda-\underline\lambda_1)t}\,\nu(d\lambda).
\]
This kernel is completely monotone since
$\lambda-\underline\lambda_1\ge0$ on $\operatorname{supp}\nu$. The map $t\mapsto \frac{e^{\underline\lambda_1t}K_1(t)}{\Lambda_1}$ belongs to $C^1(\mathbb R_+)$.
By \cite[Chapter~5, Theorems~5.2 and~5.4]{GripenbergLondenStaffans1990},
the unique locally finite signed measure satisfying the resolvent identity
for this kernel is nonnegative. Hence $J\ge0$. We now prove that $H\ge 0$.
The integrability in \eqref{eq:initial-spectral-regularity}
and the representation in \eqref{eq:initial-spectral-representation} give
\[
\forall t\ge0,\qquad
 \dot g_1(t)=-\int_{\lambda\in(0,\infty)}\lambda e^{-\lambda t}R_{1,0}(\lambda)\,\nu(d\lambda).
\]
Together with the definitions in Notations~\ref{not:return-kernel}
and~\ref{not:initial-cumulative}, this yields, for almost every $t\ge0$,
\begin{align*}
H_t
&=\beta_0\underline\lambda_1
+\beta_2\left(\frac{\dot{\bar R}_{2,t}}{2\sqrt{\bar R_{2,t}}}
+\underline\lambda_1\sqrt{\bar R_{2,t}}\right)+\beta_1\left(\dot g_1(t)+\underline\lambda_1g_1(t)\right)
-\frac{\sigma_0}{\Lambda_1}
\left(K_1'(t)+\underline\lambda_1K_1(t)\right)\\
&=\beta_0\underline\lambda_1
+\beta_2\left(\frac{\dot{\bar R}_{2,t}}{2\sqrt{\bar R_{2,t}}}
+\underline\lambda_1\sqrt{\bar R_{2,t}}\right)-\int_{\lambda\in(0,\infty)}e^{-\lambda t}
(\lambda-\underline\lambda_1)
\left(\beta_1R_{1,0}(\lambda)
-\frac{\sigma_0\lambda}{\Lambda_1}\right)\nu(d\lambda)\\
&=\beta_0\underline\lambda_1
+\beta_2\left(\frac{\dot{\bar R}_{2,t}}{2\sqrt{\bar R_{2,t}}}
+\underline\lambda_1\sqrt{\bar R_{2,t}}\right)
-\int_{\lambda\in(0,\infty)}e^{-\lambda t}\eta_0(\lambda)\,\hat\nu(d\lambda).
\end{align*}
Since~\eqref{eq:variance-lower-bound} gives $\bar R_{2,t}>0$, the chain rule, the dynamics~\eqref{eq:variance-dynamics}, the definition of $\bar R_2$ in~\eqref{PDV}, and conditions~\eqref{not:def-m} and~\eqref{eq:initial-variance-growth} in Assumption~\ref{assumption:positivity} give, for almost every $t\ge0$,
\begin{align}
\frac{\dot{\bar R}_{2,t}}{2\sqrt{\bar R_{2,t}}}
+\underline\lambda_1\sqrt{\bar R_{2,t}}
&=\frac{1}{2\sqrt{\bar R_{2,t}}}
\left(\dot g_2(t)+2\underline\lambda_1g_2(t)+K_2(t,t)\sigma_t^2 +\int_0^t
(\partial_tK_2(s,t)+2\underline\lambda_1K_2(s,t))\sigma_s^2\,ds\right)\notag\\
&\ge\frac{m}{2\sqrt{\bar R_{2,t}}}
\left(g_2(t)+\int_0^tK_2(s,t)\sigma_s^2\,ds\right)\notag\\
&=\frac{m}{2}\sqrt{\bar R_{2,t}}.\label{eq:variance-sqrt-drift-bound}
\end{align}
Taking square roots in~\eqref{eq:variance-lower-bound} gives
\[
\sqrt{\bar R_{2,t}}\ge\sqrt{\bar R_{2,0}}
e^{-(\underline\lambda_1-\frac{m}{2})t}.
\]
For $t>0$ and $\lambda>0$, we have
\[
e^{-\lambda t}=t\int_{x\in[\lambda,\infty)}e^{-xt}\,dx.
\]
The integrability in \eqref{eq:initial-d-integrability},
Fubini's theorem, and the definition of $D$ in Notation~\ref{not:initial-cumulative} yield
\begin{align}
\forall t>0,\qquad
\int_{\lambda\in(0,\infty)}e^{-\lambda t}\eta_0(\lambda)\,\hat\nu(d\lambda)
&=t\int_{x\in[\underline\lambda_1,\infty)}e^{-xt}
\left(\int_{\lambda\in(0,x]}\eta_0(\lambda)\,\hat\nu(d\lambda)\right)dx\notag\\
&\le Dt\int_{x\in[\underline\lambda_1,\infty)}e^{-xt}\,dx=De^{-\underline\lambda_1t}.\label{eq:initial-return-tail-bound}
\end{align}

\begin{itemize}
\item If $m\ge0$, for almost every $t\ge0$,
\begin{align*}
H_t
&\ge\beta_0\underline\lambda_1
+\frac{\beta_2m}{2}\sqrt{\bar R_{2,t}}-De^{-\underline\lambda_1t}\\
&\ge\beta_0\underline\lambda_1
+e^{-\underline\lambda_1t}
\left(\frac{\beta_2m}{2}\sqrt{\bar R_{2,0}}e^{\frac{mt}{2}}-D\right)\\
&\ge\beta_0\underline\lambda_1
+e^{-\underline\lambda_1t}
\left(\frac{\beta_2m}{2}\sqrt{\bar R_{2,0}}-D\right)\\
&\ge\beta_0\underline\lambda_1
-e^{-\underline\lambda_1t}
\left(D-\frac{\beta_2m}{2}\sqrt{\bar R_{2,0}}\right)_+\\
&\ge\beta_0\underline\lambda_1
-\left(D-\frac{\beta_2m}{2}\sqrt{\bar R_{2,0}}\right)_+.
\end{align*}
The first inequality follows from~\eqref{eq:variance-sqrt-drift-bound} and~\eqref{eq:initial-return-tail-bound}; the second follows from~\eqref{eq:variance-lower-bound}, since $\beta_2m\ge0$. The last three use $e^{\frac{mt}{2}}\ge1$ and $0<e^{-\underline\lambda_1t}\le1$.
\item If $m<0$, the dynamics~\eqref{eq:variance-dynamics} and the definition of $\bar R_2$ in~\eqref{PDV} give, for almost every $t\ge0$,
\[
\begin{aligned}
\dot{\bar R}_{2,t}+2\underline\lambda_1\bar R_{2,t}
&=\dot g_2(t)+2\underline\lambda_1g_2(t)+K_2(t,t)\sigma_t^2 +\int_0^t\left(\partial_tK_2(s,t)+2\underline\lambda_1K_2(s,t)\right)\sigma_s^2\,ds\\
&\ge\dot g_2(t)+2\underline\lambda_1g_2(t)\\
&\ge m g_2(t).
\end{aligned}
\]
The first inequality uses $K_2(t,t)\ge0$ from Assumption~\ref{assumption3} and the additional kernel bound in Assumption~\ref{assumption:positivity}~\myref{item:positivity-variance}; the second follows from~\eqref{eq:initial-variance-growth}.
Proposition~\ref{prop:existence-uniqueness-moments}, \eqref{eq:initial-values}, and Assumption~\ref{assumption:positivity}~\myref{item:positivity-variance} give $\bar R_{2,t}\ge g_2(t)>0$ and $g_2(t)\le g_2(0)=\bar R_{2,0}$. Since $m<0$, we obtain
\[
\frac{\dot{\bar R}_{2,t}+2\underline\lambda_1\bar R_{2,t}}{2\sqrt{\bar R_{2,t}}}
\ge\frac{m g_2(t)}{2\sqrt{\bar R_{2,t}}}
\ge\frac{m}{2}\sqrt{g_2(t)}
\ge\frac{m}{2}\sqrt{g_2(0)}
=\frac{m}{2}\sqrt{\bar R_{2,0}}.
\]
By~\eqref{eq:initial-return-tail-bound} and $D\ge0$ from Notation~\ref{not:initial-cumulative}, $\int_{\lambda\in(0,\infty)}e^{-\lambda t}\eta_0(\lambda)\,\hat\nu(d\lambda)\le De^{-\underline\lambda_1t}\le D$. Hence
$$
H_t\ge\beta_0\underline\lambda_1+\frac{\beta_2m}{2}\sqrt{\bar R_{2,0}}-D.
$$
\end{itemize}
Thus both cases give, for almost every $t\ge0$,
\begin{equation}\label{eq:positivity-H-bound}
H_t\ge\beta_0\underline\lambda_1
-\left(D-\frac{\beta_2m}{2}\sqrt{\bar R_{2,0}}\right)_+\ge0,
\end{equation}
by~\eqref{Dm} and $\beta_0\underline\lambda_1\ge0$.

\medskip

\noindent\textit{Step 2: Grönwall's lemma.}
It\^o's formula, the semimartingale dynamics~\eqref{eq:resolvent-sigma-dynamics}, and the dynamics~\eqref{eq:def-Z} of $Z$ give
\begin{align*}
d(Z_t^{-1})
&=(\lambda_{1,\mathrm{eff}}+\beta_1^2\Lambda_1^2)
Z_t^{-1}\,dt-\beta_1\Lambda_1Z_t^{-1}\,dW_t,\\
d[\sigma,Z^{-1}]_t
&=-\beta_1^2\Lambda_1^2\sigma_tZ_t^{-1}\,dt.
\end{align*}
Hence, using again \eqref{eq:resolvent-sigma-dynamics} together with \eqref{eq:resolvent-H}, we obtain
\begin{equation}\label{eq:resolvent-Zsigma}
Z_t^{-1}\sigma_t=\sigma_0+\int_0^t Z_s^{-1}
\left(((K_1'*\mathcal R)'*\sigma)(s)+H_s+(J*H)(s)\right)ds.
\end{equation}
Since $H,J\ge0$ and $\sigma_0\ge0$, Fubini's theorem and
the bounds~\eqref{eq:resolvent-memory-bound} imply, for $T>0$,
\begin{align*}
\forall t\in[0,T],\qquad
\sigma_t^-
&\le\int_{s\in[0,t]}\frac{Z_t}{Z_s}
\int_{u\in[0,s]}(K_1'*\mathcal R)'(s-u)\sigma_u^-\,du\,ds\\
&\le\left(\sup_{(u,v)\in[0,T]^2}\frac{Z_u}{Z_v}\right)
\int_{u\in[0,t]}\sigma_u^-
\left(\int_{s\in[u,t]}(K_1'*\mathcal R)'(s-u)\,ds\right)du\\
&=\left(\sup_{(u,v)\in[0,T]^2}\frac{Z_u}{Z_v}\right)
\int_{u\in[0,t]}\sigma_u^-
\left(\int_{r\in[0,t-u]}(K_1'*\mathcal R)'(r)\,dr\right)du\\
&\le\lambda_{1,\mathrm{eff}}
\left(\sup_{(u,v)\in[0,T]^2}\frac{Z_u}{Z_v}\right)
\int_0^t\sigma_u^-\,du.
\end{align*}
The supremum is finite because $Z$ is continuous and strictly positive by~\eqref{eq:def-Z}.
Grönwall's lemma gives $\sigma\ge0$. The identity~\eqref{eq:resolvent-Zsigma}
and the lower bound~\eqref{eq:positivity-H-bound} yield
\begin{align*}
\forall t\ge0,\qquad
\sigma_t&\ge Z_t\left(\sigma_0+\int_0^tZ_s^{-1}H_s\,ds\right)\\
&\ge Z_t\left(\sigma_0+
\left(\beta_0\underline\lambda_1
-\left(D-\frac{\beta_2m}{2}\sqrt{\bar R_{2,0}}\right)_+\right)
\int_0^tZ_s^{-1}\,ds\right)\ge0,
\end{align*}
which proves~\eqref{eq:positivity-pathwise-bound}.
For every $T>0$, the definition~\eqref{eq:C2T} and condition~\eqref{eq:K2-local-bounds} give
\[
C_{2,T}\le\mathcal R([0,T])
\sup_{r\in[0,T]}\sup_{u\in[0,r]}K_2(u,r)<\infty.
\]
Theorem~\ref{thm:GL_positive_moments_compact} therefore implies that
$S$ is a bounded martingale on $[0,T]$ under $\mathbb P$.
\end{proof}

\subsection{Positivity and negativity with exponential kernels}

\begin{notation}\label{not:exponential-parameters}
For the inputs~\eqref{eq:exponential-inputs}, the corresponding factors are defined, for $t\ge0$,
$1\le i\le M$ and $1\le j\le N$, by
\begin{align}\label{eq:exponential-factors}
R_{1,t}(\lambda_{1,i})&=e^{-\lambda_{1,i}t}R_{1,0}(\lambda_{1,i})
+\int_0^t\lambda_{1,i}e^{-\lambda_{1,i}(t-s)}\sigma_s\,dW_s,\notag\\
R_{2,t}(\lambda_{2,j})&\mathrel{:=}e^{-\lambda_{2,j}t}R_{2,0}(\lambda_{2,j})
+\int_0^t\lambda_{2,j}e^{-\lambda_{2,j}(t-s)}\sigma_s^2\,ds,
\end{align}
so that $\bar R_{1,t}=\sum_{i=1}^M\theta_{1,i}R_{1,t}(\lambda_{1,i})$
and $\bar R_{2,t}=\sum_{j=1}^N\theta_{2,j}R_{2,t}(\lambda_{2,j})$.
For the return kernel $K_1$, we have
\[
\nu=\sum_{i=1}^M\theta_{1,i}\delta_{\lambda_{1,i}},\qquad
\Lambda_1=\sum_{i=1}^M\theta_{1,i}\lambda_{1,i},\qquad
\underline\lambda_1=\lambda_{1,M}.
\]
By Notation~\ref{not:initial-cumulative}, $\hat\nu(\{\lambda_{1,M}\})=0$. We set
\[
\forall i\in\{1,\ldots,M-1\},\qquad
\hat\lambda_i:=\hat\nu(\{\lambda_{1,i}\})=\frac{\theta_{1,i}\lambda_{1,i}(\lambda_{1,i}-\lambda_{1,M})}{\Lambda_1}.
\]
The definitions give
\[
\begin{gathered}
\hat\nu=\sum_{i=1}^{M-1}\hat\lambda_i\delta_{\lambda_{1,i}},\qquad
\lambda_{1,\mathrm{eff}}=\lambda_{1,M}+\sum_{i=1}^{M-1}\hat\lambda_i,\qquad
D=\max_{1\le p\le M}\sum_{i=p}^{M-1}\hat\lambda_i\eta_0(\lambda_{1,i}).
\end{gathered}
\]
All sums over $1\le i\le M-1$ equal zero when $M=1$; in this case, $\hat\nu=0$ and $D=0$.
We write $\mathbf R_{2,t}\mathrel{:=}(R_{2,t}(\lambda_{2,j}))_{1\le j\le N}$ and
\[
\Lambda_2\mathrel{:=}\sum_{j=1}^N\theta_{2,j}\lambda_{2,j},\qquad
\lambda_{2,\mathrm{eff}}\mathrel{:=}
\frac{\sum_{j=1}^N\theta_{2,j}\lambda_{2,j}^2}{\Lambda_2},\qquad
\lambda_2(\mathbf r)\mathrel{:=}
\frac{\sum_{j=1}^N\theta_{2,j}\lambda_{2,j}r_j}{\sum_{j=1}^N\theta_{2,j}r_j}
\]
for $\mathbf r=(r_1,\ldots,r_N)\in\mathbb R^N$ such that
$\sum_{j=1}^N\theta_{2,j}r_j>0$.
\end{notation}

For exponential kernels, Theorem~\ref{prop:pos-MR1-generalK2-vector0} yields explicit conditions on the rates and initial factors.

\begin{prop2}\label{prop:exponential-positivity}
Suppose that $K_1$ and $K_2$ are respectively $M$- and $N$-exponential kernels, with deterministic initial factors and inputs as in Notation~\ref{not:exponential-parameters}, with $\sigma_0\ge0$. If
\begin{align}
2\lambda_{1,M}&\ge\lambda_{2,\mathrm{eff}},\label{eq:exponential-positive-rate}\\
\max_{1\le p\le M}\sum_{i=p}^{M-1}\hat\lambda_i\eta_0(\lambda_{1,i})&\le\beta_0\lambda_{1,M}
+\beta_2\left(\lambda_{1,M}-\frac{\max\{\lambda_2(\mathbf R_{2,0}),\lambda_{2,\mathrm{eff}}\}}{2}\right)\sqrt{\bar R_{2,0}},\label{eq:exponential-positive-initial}
\end{align}
then $\mathbb P(\forall t\ge0,\ \sigma_t\ge0)=1$, and $S$ is a bounded martingale on every finite interval under $\mathbb P$.
\end{prop2}
\begin{proof}
The exponential kernels and initial inputs satisfy Assumptions~\ref{assumption4} and~\ref{assumption3}, with $g_2>0$ by~\eqref{eq:exponential-input-positivity}. Proposition~\ref{prop:existence-uniqueness-moments} therefore gives a pathwise unique continuous solution of~\eqref{PDV}. Definition~\ref{NexponentialKernel} gives $\theta_{2,j},\lambda_{2,j}>0$ and $R_{2,0}(\lambda_{2,j})\ge0$ for every $j\in\{1,\ldots,N\}$. Together with~\eqref{eq:exponential-inputs}, this gives, for every $t\ge0$,
\begin{gather*}
K_2(0,t)\partial_t^2K_2(0,t)-(\partial_tK_2(0,t))^2=\frac{1}{2}\sum_{j,k=1}^N\theta_{2,j}\theta_{2,k}\lambda_{2,j}\lambda_{2,k}(\lambda_{2,j}-\lambda_{2,k})^2e^{-(\lambda_{2,j}+\lambda_{2,k})t}\ge0,\\
g_2(t)\ddot g_2(t)-\dot g_2(t)^2=\frac{1}{2}\sum_{j,k=1}^N\theta_{2,j}\theta_{2,k}R_{2,0}(\lambda_{2,j})R_{2,0}(\lambda_{2,k})(\lambda_{2,j}-\lambda_{2,k})^2e^{-(\lambda_{2,j}+\lambda_{2,k})t}\ge0.
\end{gather*}
Since $K_2(0,t)>0$ by Definition~\ref{NexponentialKernel} and $g_2(t)>0$ by~\eqref{eq:exponential-input-positivity}, the following quantities are well defined and satisfy
\[
\forall t\ge0,\qquad
\frac{d^2}{dt^2}\log K_2(0,t)\ge0,\qquad
\frac{d^2}{dt^2}\log g_2(t)\ge0.
\]
Thus $K_2(0,\cdot)$ and $g_2$ are log-convex. Their logarithmic derivatives are nondecreasing, so evaluation at zero gives
\begin{equation}\label{eq:exponential-log-derivative-bounds}
\begin{aligned}
\forall t\ge0,\qquad -\frac{\partial_tK_2(0,t)}{K_2(0,t)}
&\le-\frac{\partial_tK_2(0,0)}{K_2(0,0)}=\lambda_{2,\mathrm{eff}},\\
\forall t\ge0,\qquad -\frac{\dot g_2(t)}{g_2(t)}
&\le-\frac{\dot g_2(0)}{g_2(0)}=\lambda_2(\mathbf R_{2,0}).
\end{aligned}
\end{equation}
The logarithmic-derivative bounds~\eqref{eq:exponential-log-derivative-bounds} are attained at zero. Hence the constant in Assumption~\ref{assumption:positivity} is $m=2\lambda_{1,M}-\max\{\lambda_2(\mathbf R_{2,0}),\lambda_{2,\mathrm{eff}}\}$.
\begin{itemize}
\item If $m\ge0$, the first alternative in Assumption~\ref{assumption:positivity} holds.
\item If $m<0$,~\eqref{eq:exponential-positive-rate} and the first bound in~\eqref{eq:exponential-log-derivative-bounds}, using $K_2(s,t)=K_2(0,t-s)$, give, for every $0\le s\le t$,
$\partial_tK_2(s,t)+2\lambda_{1,M}K_2(s,t)\ge0$.
Definition~\ref{NexponentialKernel} and~\eqref{eq:exponential-inputs} give $g_2(t)\le g_2(0)$.
\end{itemize}
Finally,~\eqref{eq:exponential-positive-initial} is equivalent to
$D\le\beta_0\lambda_{1,M}+\frac{\beta_2m}{2}\sqrt{\bar R_{2,0}}$.
Thus Assumption~\ref{assumption:positivity} holds, and Theorem~\ref{prop:pos-MR1-generalK2-vector0} applies.

\end{proof}

\begin{rk}\label{rk:exponential-initial-rate}
For an $N$-exponential kernel $K_2$, suppose that, for some $\Delta\in(0,\infty]$ and a deterministic measurable history $(\sigma_s)_{-\Delta<s<0}$,
\[
\forall j\in\{1,\ldots,N\},\qquad
R_{2,0}(\lambda_{2,j})=\lambda_{2,j}\int_{-\Delta}^0e^{\lambda_{2,j}s}\sigma_s^2\,ds,
\qquad
\int_{-\Delta}^0e^{\lambda_{2,N}s}\sigma_s^2\,ds\in(0,\infty).
\]These assumptions give $R_{2,0}(\lambda_{2,j})\in(0,\infty)$ for every $j\in\{1,\ldots,N\}$, hence $g_2(0)>0$ by~\eqref{eq:exponential-inputs}. Thus,
\begin{equation}\label{eq:exponential-initial-rate}
\lambda_2(\mathbf R_{2,0})\le\lambda_{2,\mathrm{eff}}.
\end{equation}
Consequently, condition~\eqref{eq:exponential-positive-initial} simplifies to
\[
\max_{1\le p\le M}\sum_{i=p}^{M-1}\hat\lambda_i\eta_0(\lambda_{1,i})
\le\beta_0\lambda_{1,M}
+\beta_2\left(\lambda_{1,M}-\frac{\lambda_{2,\mathrm{eff}}}{2}\right)\sqrt{\bar R_{2,0}},
\]
since $\max\{\lambda_2(\mathbf R_{2,0}),\lambda_{2,\mathrm{eff}}\}=\lambda_{2,\mathrm{eff}}$. By~\eqref{eq:exponential-inputs}, $g_2(t)=\int_{-\Delta}^0K_2(0,t-s)\sigma_s^2\,ds$. Differentiating at zero and applying the first bound in~\eqref{eq:exponential-log-derivative-bounds} at $t=-s\ge0$ give
\[
\begin{aligned}\lambda_2(\mathbf R_{2,0})&=-\frac{\dot g_2(0)}{g_2(0)}
=\frac{\int_{-\Delta}^0-\partial_tK_2(0,-s)\sigma_s^2\,ds}{\int_{-\Delta}^0K_2(0,-s)\sigma_s^2\,ds}
\le\lambda_{2,\mathrm{eff}}\frac{\int_{-\Delta}^0K_2(0,-s)\sigma_s^2\,ds}{\int_{-\Delta}^0K_2(0,-s)\sigma_s^2\,ds}=\lambda_{2,\mathrm{eff}}.\end{aligned}
\]
\end{rk}
The next result gives conditions under which volatility becomes negative with positive probability, and hence necessary conditions for nonnegative volatility.

\begin{thm2}\label{gonegative-N-eff2}
Suppose that $K_1$ and $K_2$ are respectively $M$- and $N$-exponential kernels. We use Notations~\ref{not:initial-cumulative} and~\ref{not:exponential-parameters}. Assume that the deterministic initial factors and inputs are as in~\eqref{eq:exponential-inputs}, with $\sigma_0>0$. If at least one of the following
conditions holds:
\begin{align}
\lambda_{1,\mathrm{eff}}
&<\frac{\lambda_{2,\mathrm{eff}}}{2},
\label{eq:negativity-effective-rate}\\
\sum_{i=1}^{M-1}\hat\lambda_i\eta_0(\lambda_{1,i})
&>\beta_0\lambda_{1,M}
+\beta_2\left(
\lambda_{1,M}-\frac{\lambda_2(\mathbf R_{2,0})}{2}
\right)\sqrt{\bar R_{2,0}}.
\label{eq:negativity-initial-drift}
\end{align}
then we have
\[
\forall t>0,\qquad \mathbb P(\sigma_t<0)>0.
\]
Consequently, for every $t_0>0$,
$(S_t)_{t\in[0,t_0]}$ is a strict local martingale under
$\mathbb P$.

\end{thm2}

For $M=1$, $\hat\nu=0$, $D=0$, and $\lambda_{1,\mathrm{eff}}=\lambda_{1,1}$. Combining Proposition~\ref{prop:exponential-positivity}, Theorem~\ref{gonegative-N-eff2}, and Remark~\ref{rk:exponential-initial-rate} yields the following corollary.

\begin{cor2}\label{cor:M1-equivalence}
Suppose $M=1$ and $K_2$ is an $N$-exponential kernel, with deterministic initial factors and inputs as in Notation~\ref{not:exponential-parameters}, with $\sigma_0>0$.
\begin{enumerate}[(i)]
\myitem[(i)]{item:M1-general} The following equivalences hold:
\begin{align*}
\mathbb P(\forall t\ge0,\ \sigma_t\ge0)=1
&\quad\Longleftrightarrow\quad
\begin{cases}
2\lambda_{1,1}\ge\lambda_{2,\mathrm{eff}},\\
0\le\beta_0\lambda_{1,1}
+\beta_2\left(\lambda_{1,1}-\frac{\lambda_2(\mathbf R_{2,0})}{2}\right)\sqrt{\bar R_{2,0}}
\end{cases}\\
&\quad\Longleftrightarrow\quad
\text{$S$ is a martingale on every finite interval under $\mathbb P$.}
\end{align*}
\myitem[(ii)]{item:M1-past}
If $g_2$ is initialized from a deterministic past as in Remark~\ref{rk:exponential-initial-rate}, then
\begin{align*}
\mathbb P(\forall t\ge0,\ \sigma_t\ge0)=1
&\quad\Longleftrightarrow\quad 2\lambda_{1,1}\ge\lambda_{2,\mathrm{eff}}\\
&\quad\Longleftrightarrow\quad
\text{$S$ is a bounded martingale on every finite interval under $\mathbb P$.}
\end{align*}
\end{enumerate}

\end{cor2}

\begin{rk}
When $g_2$ is initialized from a deterministic past as in Remark~\ref{rk:exponential-initial-rate}, nonnegative volatility for $M=1$ is characterized by $2\lambda_{1,1}\ge\lambda_{2,\mathrm{eff}}$, with no further restriction on $R_{1,0}(\lambda_{1,1})$ beyond $\sigma_0>0$. For every $M\ge2$, under the initial-state hypotheses of Theorem~\ref{gonegative-N-eff2}, nonnegative volatility also requires the following restriction on the initial return factors (the reverse of~\eqref{eq:negativity-initial-drift}):
\[
\sum_{i=1}^{M-1}\hat\lambda_i\eta_0(\lambda_{1,i})
\le\beta_0\lambda_{1,M}
+\beta_2\left(\lambda_{1,M}-\frac{\lambda_2(\mathbf R_{2,0})}{2}\right)\sqrt{\bar R_{2,0}}.
\]
When the initial factors are chosen freely, this necessary condition depends on $(R_{1,0}(\lambda_{1,i}))_{1\le i\le M}$, even when $\sigma_0>0$ and $\mathbf R_{2,0}$ are fixed. Indeed, for $a>0$, the changes
\[
\begin{aligned}
R_{1,0}(\lambda_{1,1})&\mapsto R_{1,0}(\lambda_{1,1})-a,\\
R_{1,0}(\lambda_{1,M})&\mapsto R_{1,0}(\lambda_{1,M})+\frac{\theta_{1,1}}{\theta_{1,M}}a,
\end{aligned}
\]
with all other factors fixed, preserve $\bar R_{1,0}$ and hence $\sigma_0$. They leave the right-hand side unchanged and increase the left-hand side by
\[
\hat\lambda_1\frac{\beta_1\Lambda_1}{\lambda_{1,1}}(-a)
=-\beta_1\theta_{1,1}(\lambda_{1,1}-\lambda_{1,M})a
\xrightarrow[a\to+\infty]{}+\infty.
\]
\end{rk}

Table~\ref{tab:suff-necessary-conditions} below summarizes the sufficient and necessary
conditions for nonnegative volatility given by
Proposition~\ref{prop:exponential-positivity} and Theorem~\ref{gonegative-N-eff2}, together with Corollary~\ref{cor:M1-equivalence},
when $K_1$ and $K_2$ are $M$- and $N$-exponential kernels, respectively, using
Notations~\ref{not:exponential-parameters} and~\ref{not:initial-cumulative}.
\begin{table}[H]
\centering
\small
\renewcommand{\arraystretch}{1.4}
\setlength{\tabcolsep}{5pt}
\begin{tabularx}{\textwidth}{@{}l>{\centering\arraybackslash}X>{\centering\arraybackslash}X@{}}
\toprule
$K_1$ & Sufficient conditions for positivity of volatility & Necessary conditions for positivity of volatility \\
\midrule
$1$-exp & $\begin{aligned}
&2\lambda_{1,1}\ge\lambda_{2,\mathrm{eff}},\\
&0\le\beta_0\lambda_{1,1}\\
&\quad+\beta_2\left(\lambda_{1,1}-\frac{\lambda_2(\mathbf R_{2,0})}{2}\right)\sqrt{\bar R_{2,0}}
\end{aligned}$
& $\begin{aligned}
&2\lambda_{1,1}\ge\lambda_{2,\mathrm{eff}},\\
&0\le\beta_0\lambda_{1,1}\\
&\quad+\beta_2\left(\lambda_{1,1}-\frac{\lambda_2(\mathbf R_{2,0})}{2}\right)\sqrt{\bar R_{2,0}}
\end{aligned}$ \\
\midrule
$2$-exp
& $\begin{aligned}
&2\lambda_{1,2}\ge\lambda_{2,\mathrm{eff}},\\
&\max\{0,\hat\lambda_1\eta_0(\lambda_{1,1})\}\le\beta_0\lambda_{1,2}\\
&\quad+\beta_2\left(\lambda_{1,2}-\frac{\lambda_2(\mathbf R_{2,0})}{2}\right)\sqrt{\bar R_{2,0}},\\
&\max\{0,\hat\lambda_1\eta_0(\lambda_{1,1})\}\le\beta_0\lambda_{1,2}\\
&\quad+\beta_2\left(\lambda_{1,2}-\frac{\lambda_{2,\mathrm{eff}}}{2}\right)\sqrt{\bar R_{2,0}}
\end{aligned}$
& $\begin{aligned}
&\lambda_{1,\mathrm{eff}}\ge\frac{\lambda_{2,\mathrm{eff}}}{2},\\
&\hat\lambda_1\eta_0(\lambda_{1,1})\le\beta_0\lambda_{1,2}\\
&{}+\beta_2\left(\lambda_{1,2}-\frac{1}{2}\lambda_2(\mathbf R_{2,0})\right)\sqrt{\bar R_{2,0}}
\end{aligned}$ \\
\midrule
$M$-exp
& $\begin{aligned}
&2\lambda_{1,M}\ge\lambda_{2,\mathrm{eff}},\\
&\max_{1\le p\le M}\sum_{i=p}^{M-1}\hat\lambda_i\eta_0(\lambda_{1,i})\le\beta_0\lambda_{1,M}\\
&\quad+\beta_2\left(\lambda_{1,M}-\frac{\lambda_2(\mathbf R_{2,0})}{2}\right)\sqrt{\bar R_{2,0}},\\
&\max_{1\le p\le M}\sum_{i=p}^{M-1}\hat\lambda_i\eta_0(\lambda_{1,i})\le\beta_0\lambda_{1,M}\\
&\quad+\beta_2\left(\lambda_{1,M}-\frac{\lambda_{2,\mathrm{eff}}}{2}\right)\sqrt{\bar R_{2,0}}
\end{aligned}$
& $\begin{aligned}
&\lambda_{1,\mathrm{eff}}\ge\frac{\lambda_{2,\mathrm{eff}}}{2},\\
&\sum_{i=1}^{M-1}\hat\lambda_i\eta_0(\lambda_{1,i})\le\beta_0\lambda_{1,M}\\
&{}+\beta_2\left(\lambda_{1,M}-\frac{1}{2}\lambda_2(\mathbf R_{2,0})\right)\sqrt{\bar R_{2,0}}
\end{aligned}$ \\
\bottomrule
\end{tabularx}
\caption{Conditions from Proposition~\ref{prop:exponential-positivity}
and Theorem~\ref{gonegative-N-eff2} when $K_2$ is an $N$-exponential kernel.
The deterministic initial factors and inputs satisfy
\eqref{eq:exponential-inputs} and $\bar R_{2,0}=g_2(0)>0$ as in~\eqref{eq:initial-values}.
The sufficient conditions assume $\sigma_0\ge0$, and the necessary
conditions assume $\sigma_0>0$.
The sufficient and necessary conditions coincide for $M=1$ when $\sigma_0>0$.
For $M\ge2$, the sufficient and necessary rate conditions involve
$\lambda_{1,M}$ and $\lambda_{1,\mathrm{eff}}$, respectively;
the sufficient conditions on the initial factors involve both
$\lambda_{2,\mathrm{eff}}$ and $\lambda_2(\mathbf R_{2,0})$, while the necessary condition involves $\lambda_2(\mathbf R_{2,0})$.
For $M\ge3$, the sufficient conditions bound the maximum of zero and all tail sums
$\sum_{i=p}^{M-1}\hat\lambda_i\eta_0(\lambda_{1,i})$,
whereas the necessary condition bounds only the full sum corresponding to $p=1$.}
\label{tab:suff-necessary-conditions}
\end{table}

\begin{rk}
When the Guyon--Lekeufack model~\eqref{PDV} is initialized from observed past returns, its initial squared-return factors are determined by the historical price path rather than chosen freely; see \cite[Remark~2.2]{GazzaniGuyon2025}. For $n$ daily returns, we set $\Delta\mathrel{:=}\frac{n}{252}$ years for the historical window. From the prices $S_{-\Delta},S_{-\Delta+\frac{1}{252}},\ldots,S_0$, we define the returns by
\[
\forall k\in\{1,\ldots,n\},\qquad
r_k\mathrel{:=}\frac{S_{\frac{k-n}{252}}}{S_{\frac{k-1-n}{252}}}-1,
\]
and we define the piecewise constant variance history by
\[
\forall s\in[-\Delta,0),\qquad
\sigma_s^2\mathrel{:=}252\sum_{k=1}^n r_k^2
\mathbf{1}_{\left[\frac{k-1-n}{252},\frac{k-n}{252}\right)}(s).
\]
If $\sum_{k=1}^n r_k^2>0$, this history satisfies the integrability assumption in Remark~\ref{rk:exponential-initial-rate} and yields $R_{2,0}(\lambda_{2,j})=\lambda_{2,j}\int_{-\Delta}^0e^{\lambda_{2,j}s}\sigma_s^2\,ds$, hence~\eqref{eq:exponential-initial-rate}.
\end{rk}

\begin{example}\label{ex:market-bound}\label{ex:market-positivity}
We use the following parameters and initial factors, with the parameters based on the market-smile calibration in \cite[Table~8]{GuyonLekeufack2023}:
\[
\begin{aligned}
(\beta_0,\beta_1,\beta_2)&=(0.026,-0.138,0.69),\\
(\lambda_{1,1},\lambda_{1,2},\theta_{1,1},\theta_{1,2})&=(62.11,32.25,0.77,0.23),\\
(\lambda_{2,1},\lambda_{2,2},\theta_{2,1},\theta_{2,2})&=(9.57,3.51,0.01,0.99),\\
\left(R_{1,0}(\lambda_{1,1}),R_{1,0}(\lambda_{1,2})\right)&=(0.2988,0.2397),\\
\left(R_{2,0}(\lambda_{2,1}),R_{2,0}(\lambda_{2,2})\right)&=(0.016,0.02).
\end{aligned}
\]
These parameters and initial factors, together with~\eqref{eq:initial-values}, give
\[
\bar R_{2,0}=0.01996>0,\qquad \sigma_0\simeq0.08412>0,\qquad
\lambda_2(\mathbf R_{2,0})\simeq3.55858
<\lambda_{2,\mathrm{eff}}\simeq3.67242.
\]
Moreover, $\eta_0(\lambda_{1,1})\simeq-0.12080<0$ implies $D=0$ by
Notation~\ref{not:exponential-parameters}.
The sufficient conditions in the second row of Table~\ref{tab:suff-necessary-conditions}
are satisfied, since
\[
\begin{aligned}
\lambda_{1,2}=32.25&>\frac{\lambda_{2,\mathrm{eff}}}{2}\simeq1.83621,\\
\max\{0,\hat\lambda_1\eta_0(\lambda_{1,1})\}=0
&<\beta_0\lambda_{1,2}
+\beta_2\left(\lambda_{1,2}-\frac{\lambda_{2,\mathrm{eff}}}{2}\right)
\sqrt{\bar R_{2,0}}\simeq3.80333.
\end{aligned}
\]
Thus Theorem~\ref{prop:pos-MR1-generalK2-vector0} applies with
$m=2\lambda_{1,2}-\lambda_{2,\mathrm{eff}}\simeq60.82758>0$
and gives $\sigma_t>0$ for every $t\ge0$, almost surely. With $S_0=1$ and $T=\frac{16}{252}$, bound~\eqref{eq:GL_asset_pathwise_upper_bound} gives $B_T\simeq1.40630$.
With $2\times10^6$ Monte Carlo paths and a time step of $\frac{1}{25200}$ years, the call-price estimate first vanished at the tested strike $K=1.05$: no simulated terminal value exceeded this strike.
\end{example}

For the proof of Theorem~\ref{gonegative-N-eff2}, $K_1$ and $K_2$ are respectively $M$- and $N$-exponential kernels as in Definition~\ref{NexponentialKernel}, with $M,N\ge1$ and $\sigma_0>0$.
We define the processes
\begin{equation}\label{eq:negativity-z-definition}
\forall t\ge0,\quad\forall i\in\{1,\ldots,M-1\},\qquad
z_{i,t}:=\beta_1\left(
\frac{\Lambda_1}{\lambda_{1,i}}R_{1,t}(\lambda_{1,i})-\bar R_{1,t}\right).
\end{equation}
By~\eqref{eq:exponential-factors} and~\eqref{eq:negativity-z-definition}, they have finite variation, since their Brownian coefficients are
\[
\beta_1\left(
\frac{\Lambda_1}{\lambda_{1,i}}\lambda_{1,i}-\Lambda_1
\right)\sigma_t=0.
\]
We write $z_t\mathrel{:=}(z_{i,t})_{1\le i\le M-1}$.
By~\eqref{eq:exponential-factors} and Notation~\ref{not:exponential-parameters}, the drift of $\beta_1\bar R_{1,t}$ is
\begin{align}
-\beta_1\sum_{i=1}^M\theta_{1,i}\lambda_{1,i}R_{1,t}(\lambda_{1,i})
&=-\lambda_{1,M}\beta_1\bar R_{1,t}
-\sum_{i=1}^{M-1}\hat\lambda_i
\left(z_{i,t}+\beta_1\bar R_{1,t}\right)\notag\\
&=-\lambda_{1,\mathrm{eff}}
(\sigma_t-\beta_0-\beta_2\sqrt{\bar R_{2,t}})
-\sum_{i=1}^{M-1}\hat\lambda_iz_{i,t}.
\label{eq:negativity-return-drift}
\end{align}
Together with the factor dynamics~\eqref{eq:exponential-factors} and $\bar R_{2,t}\ge g_2(t)>0$ from~\eqref{PDV} and~\eqref{eq:exponential-input-positivity}, the chain rule gives
\begin{subequations}\label{eq:negativity-lift}
\begin{align}
d\sigma_t
&=\mu_\sigma(\sigma_t,z_t,\mathbf R_{2,t})\,dt
+\beta_1\Lambda_1\sigma_t\,dW_t,\label{eq:negativity-sigma-dynamics}\\
\forall i\in\{1,\ldots,M-1\},\qquad\dot z_{i,t}
&=-\lambda_{1,i}z_{i,t}+\sum_{k=1}^{M-1}\hat\lambda_kz_{k,t}
+(\lambda_{1,\mathrm{eff}}-\lambda_{1,i})
(\sigma_t-\beta_0-\beta_2\sqrt{\bar R_{2,t}}),\label{eq:negativity-z-dynamics}\\
\forall j\in\{1,\ldots,N\},\qquad\dot R_{2,t}(\lambda_{2,j})
&=\lambda_{2,j}\left(\sigma_t^2-R_{2,t}(\lambda_{2,j})\right),\label{eq:negativity-variance-dynamics}
\end{align}
\end{subequations}
where we define the drift coefficient by
\begin{align}\label{eq:negativity-drift-definition}
\mu_\sigma(\sigma,z,\mathbf r)
&\mathrel{:=}-\lambda_{1,\mathrm{eff}}\sigma
+\lambda_{1,\mathrm{eff}}\beta_0-\sum_{i=1}^{M-1}\hat\lambda_iz_i
+\beta_2\left(\lambda_{1,\mathrm{eff}}-
\frac{\lambda_2(\mathbf r)}{2}\right)\sqrt{\sum_{j=1}^N\theta_{2,j}r_j}
+\frac{\beta_2\Lambda_2\sigma^2}{2\sqrt{\sum_{j=1}^N\theta_{2,j}r_j}}.
\end{align}
This coefficient is well defined and smooth for $\sigma\in\mathbb R$,
$z\in\mathbb R^{M-1}$, and $\mathbf r\in\mathbb R^N$ with
$\sum_{j=1}^N\theta_{2,j}r_j>0$. For $i\in\{1,\ldots,M-1\}$, the drift in~\eqref{eq:negativity-z-dynamics} follows from~\eqref{eq:negativity-z-definition} and~\eqref{eq:negativity-return-drift}:
\begin{align*}
\dot z_{i,t}
&=-\beta_1\Lambda_1R_{1,t}(\lambda_{1,i})
+\beta_1\sum_{k=1}^M\theta_{1,k}\lambda_{1,k}R_{1,t}(\lambda_{1,k})\\
&=-\lambda_{1,i}z_{i,t}
+(\lambda_{1,\mathrm{eff}}-\lambda_{1,i})\beta_1\bar R_{1,t}
+\sum_{k=1}^{M-1}\hat\lambda_kz_{k,t},
\end{align*}
using $\beta_1\bar R_{1,t}=\sigma_t-\beta_0-\beta_2\sqrt{\bar R_{2,t}}$. When $M=1$, $z$ and all sums over $1\le i\le M-1$ are omitted. For $T>0$ and a deterministic control $u\in L^2([0,T])$,
we denote by $(\sigma^u,z^u,\mathbf R_2^u)$ the solution of
\eqref{eq:negativity-lift} on its interval of existence, with $dW_t$
replaced by $u_t\,dt$ and the same initial state. We write
\[
\mathbf R_{2,t}^u\mathrel{:=}(R_{2,t}^u(\lambda_{2,j}))_{1\le j\le N},\qquad
\bar R_{2,t}^u\mathrel{:=}\sum_{j=1}^N\theta_{2,j}R_{2,t}^u(\lambda_{2,j}).
\]
For $h>0$ and a continuous positive path $s^h$ on $[0,h]$
with $s_0^h=\sigma_0$, we define the functions $\widetilde z$ and
$\widetilde{\mathbf R}_2\mathrel{:=}(\widetilde R_2(\lambda_{2,j}))_{1\le j\le N}$ as the solutions of
\eqref{eq:negativity-z-dynamics} and
\eqref{eq:negativity-variance-dynamics}
with $\sigma=s^h$ and initial values $z_0$ and $\mathbf R_{2,0}$.

\begin{rk}
In~\eqref{eq:negativity-z-dynamics}, the Brownian terms cancel and the drift contains $\sigma$; in~\eqref{eq:negativity-variance-dynamics}, the variance-factor drift contains $\sigma^2$. If a volatility path has size
$h^{-\frac{1}{2}}$ on $[0,h]$, these terms give, up to constants
independent of $h$,
\begin{align*}
|\widetilde z_{i,h}-z_{i,0}|
&\approx h\times h^{-\frac{1}{2}}=\sqrt h,\qquad
\widetilde R_{2,h}(\lambda_{2,j})-R_{2,0}(\lambda_{2,j})
\approx h\times(h^{-\frac{1}{2}})^2=1.
\end{align*}
When $\widetilde z$ and $\widetilde{\mathbf R}_2$ stay bounded, the terms in their
drifts that do not contain $\sigma$ or $\sigma^2$ contribute $O(h)$.
Thus the variance factors can change by an amount of order one while
$\widetilde z$ remains close to $z_0$. The proof uses a positive path $s^h$ on $[0,h]$ with $s_0^h=\sigma_0$
and
\[
s_h^h=h^2,\qquad
\int_0^h s_t^h\,dt\xrightarrow[h\downarrow0]{}0,\qquad
\int_0^h(s_t^h)^2\,dt\xrightarrow[h\downarrow0]{}A
\]
for a chosen $A\ge0$. The controlled factors then satisfy
\[
\widetilde z_h\xrightarrow[h\downarrow0]{}z_0,\qquad
\widetilde R_{2,h}(\lambda_{2,j})
\xrightarrow[h\downarrow0]{} R_{2,0}(\lambda_{2,j})+\lambda_{2,j}A.
\]
We choose $A$ so that
\[
\mu_\sigma\left(0,z_0,
(R_{2,0}(\lambda_{2,j})+\lambda_{2,j}A)_{1\le j\le N}\right)<0.
\]
Since $s_h^h=h^2=o(h)$, setting $u=0$ on $[h,2h]$ makes $\sigma^u$
negative for all sufficiently small $h$. A control on $[2h,T]$
keeps $\sigma^u$ at its negative value at time $2h$.
\end{rk}

\begin{proof}[Proof of Theorem~\ref{gonegative-N-eff2}]
We fix $T>0$.
\medskip

\noindent\textit{Step 1: Construction of a control.}
Since $z_{i,0}=\eta_0(\lambda_{1,i})+\beta_0+\beta_2\sqrt{\bar R_{2,0}}$ (by \eqref{eq:negativity-z-definition} and Notation~\ref{not:initial-cumulative}) and
$\lambda_{1,\mathrm{eff}}-\sum_{i=1}^{M-1}\hat\lambda_i=\lambda_{1,M}$, we have
\begin{align}
\mu_\sigma(0,z_0,\mathbf R_{2,0})
&=\beta_0\left(\lambda_{1,\mathrm{eff}}-\sum_{i=1}^{M-1}\hat\lambda_i\right)
-\sum_{i=1}^{M-1}\hat\lambda_i\eta_0(\lambda_{1,i})\notag\\
&\quad+\beta_2\left(\lambda_{1,\mathrm{eff}}-\sum_{i=1}^{M-1}\hat\lambda_i
-\frac{\lambda_2(\mathbf R_{2,0})}{2}\right)\sqrt{\bar R_{2,0}}\notag\\
&=\beta_0\lambda_{1,M}-\sum_{i=1}^{M-1}\hat\lambda_i\eta_0(\lambda_{1,i})
+\beta_2\left(\lambda_{1,M}-
\frac{\lambda_2(\mathbf R_{2,0})}{2}\right)\sqrt{\bar R_{2,0}}.
\label{eq:negativity-boundary-drift}
\end{align}
Moreover, we have
\begin{align*}
\sum_{j=1}^N\theta_{2,j}
\left(R_{2,0}(\lambda_{2,j})+\lambda_{2,j}A\right)
&=\bar R_{2,0}+\Lambda_2A,\\
\lambda_2\left((R_{2,0}(\lambda_{2,j})+\lambda_{2,j}A)_{1\le j\le N}\right)
&=\frac{\lambda_2(\mathbf R_{2,0})\bar R_{2,0}
+\Lambda_2\lambda_{2,\mathrm{eff}}A}{\bar R_{2,0}+\Lambda_2A}
\xrightarrow[A\to\infty]{}\lambda_{2,\mathrm{eff}}.
\end{align*}
Substitution in \eqref{eq:negativity-drift-definition} gives, for $A>0$,

\begin{equation}\label{eq:negativity-drift-limit}
\begin{aligned}
\frac{\mu_\sigma\left(0,z_0,
(R_{2,0}(\lambda_{2,j})+\lambda_{2,j}A)_{1\le j\le N}\right)}{\sqrt A}
&=\frac{\lambda_{1,\mathrm{eff}}\beta_0
-\sum_{i=1}^{M-1}\hat\lambda_iz_{i,0}}{\sqrt A}\\
&\quad+\beta_2\left(\lambda_{1,\mathrm{eff}}
-\frac{\lambda_2(\mathbf R_{2,0})\bar R_{2,0}
+\Lambda_2\lambda_{2,\mathrm{eff}}A}{2(\bar R_{2,0}+\Lambda_2A)}\right)
\sqrt{\Lambda_2+\frac{\bar R_{2,0}}{A}}\\
&\xrightarrow[A\to\infty]{}
\beta_2\sqrt{\Lambda_2}
\left(\lambda_{1,\mathrm{eff}}-
\frac{\lambda_{2,\mathrm{eff}}}{2}\right).
\end{aligned}
\end{equation}
If condition~\eqref{eq:negativity-initial-drift} of Theorem~\ref{gonegative-N-eff2} holds, $A=0$ makes the drift
negative by the identity~\eqref{eq:negativity-boundary-drift}.
If condition~\eqref{eq:negativity-effective-rate} holds, every sufficiently
large $A$ makes the drift negative by the limit~\eqref{eq:negativity-drift-limit}.
Thus there exists $A\ge0$ such that
\begin{equation}\label{eq:negativity-target-drift}
\mu_\sigma\left(0,z_0,
(R_{2,0}(\lambda_{2,j})+\lambda_{2,j}A)_{1\le j\le N}\right)<0.
\end{equation}
For $h\in(0,1]$, we define
\[
\forall t\in[0,h],\qquad
s_t^h\mathrel{:=}\sigma_0\left(1-\frac{t}{h}\right)+ht+\sqrt{\frac{2A}{h}}\sin\left(\frac{\pi t}{h}\right).
\]
For fixed $h$ and every $t\in[0,h]$,
\[
s_h^h=h^2,\qquad
s_t^h\ge \sigma_0\left(1-\frac{t}{h}\right)+ht
\ge \min\{\sigma_0,h^2\}>0.
\]
By~\eqref{eq:negativity-variance-dynamics} and~\eqref{eq:initial-values},
\[
\forall t\in[0,h],\qquad
\sum_{j=1}^N\theta_{2,j}\widetilde R_{2,t}(\lambda_{2,j})\ge e^{-\lambda_{2,1}h}\bar R_{2,0}>0.
\]
Thus $\mu_\sigma(s_t^h,\widetilde z_t,\widetilde{\mathbf R}_{2,t})$ is well defined and $\beta_1\Lambda_1s_t^h\ne0$. The control on $[0,h]$ is therefore well defined by
\[
\forall t\in[0,h],\qquad
u_t\mathrel{:=}\frac{\dot s_t^h-\mu_\sigma(s_t^h,\widetilde z_t,\widetilde{\mathbf R}_{2,t})}{\beta_1\Lambda_1s_t^h}.
\]
It is continuous, hence bounded on $[0,h]$.
This gives $\sigma_t^u=s_t^h$ on $[0,h]$. Uniqueness of the factor equations then gives $z_t^u=\widetilde z_t$ and $\mathbf R_{2,t}^u=\widetilde{\mathbf R}_{2,t}$ for $0\le t\le h$.
Moreover,
\begin{align*}
\int_0^h s_t^h\,dt
&=\frac{\sigma_0h}{2}+\frac{h^3}{2}
+\frac{2\sqrt{2A}}{\pi}\sqrt h
=O(\sqrt h),\\
\int_0^h(s_t^h)^2\,dt
&=A+\frac{2\sqrt{2A}}{\pi}(\sigma_0+h^2)\sqrt h
+\frac{h}{3}(\sigma_0^2+\sigma_0h^2+h^4)=A+O(\sqrt h).
\end{align*}
Hence, the variance factors satisfy
\begin{align}
\widetilde R_{2,h}(\lambda_{2,j})
&=e^{-\lambda_{2,j}h}R_{2,0}(\lambda_{2,j})
+\lambda_{2,j}\int_0^h
e^{-\lambda_{2,j}(h-t)}(s_t^h)^2\,dt
\notag\\
&=R_{2,0}(\lambda_{2,j})
+\lambda_{2,j}\int_0^h(s_t^h)^2\,dt+O(h)
\notag\\
&=R_{2,0}(\lambda_{2,j})+\lambda_{2,j}A+O(\sqrt h).
\label{eq:negativity-variance-shift}
\end{align}
The $O(h)$ remainder in \eqref{eq:negativity-variance-shift} follows
from $0\le1-e^{-x}\le x$ for $x\ge0$ and
\begin{align*}
\left|\widetilde R_{2,h}(\lambda_{2,j})-R_{2,0}(\lambda_{2,j})
-\lambda_{2,j}\int_0^h(s_t^h)^2\,dt\right|\le\lambda_{2,j}hR_{2,0}(\lambda_{2,j})
+\lambda_{2,j}^2h\int_0^h(s_t^h)^2\,dt
=O(h).
\end{align*}
Moreover, we have
\[
\forall t\in[0,h],\ \forall j\in\{1,\ldots,N\},\qquad
0\le \widetilde R_{2,t}(\lambda_{2,j})
\le R_{2,0}(\lambda_{2,j})
+\lambda_{2,j}\int_0^h(s_r^h)^2\,dr\le C.
\]
In the estimates on $[0,h]$, $C$ denotes a constant depending only on
$A$, the parameters and the initial state.
The equation~\eqref{eq:negativity-z-dynamics} for $z$ therefore gives
\[
|\dot{\widetilde z}_t|\le C(1+|\widetilde z_t|+s_t^h).
\]
Since $\int_0^h s_t^h\,dt=O(\sqrt h)$, Grönwall's lemma yields
\[
\sup_{t\in[0,h]}|\widetilde z_t|\le C,
\qquad
|\widetilde z_h-z_0|
\le C\int_0^h(1+|\widetilde z_t|+s_t^h)\,dt
=O(\sqrt h).
\]
Since $\bar R_{2,0}+\Lambda_2A>0$ by~\eqref{eq:initial-values}, \eqref{eq:negativity-drift-definition} shows that the coefficients in~\eqref{eq:negativity-lift} are bounded and Lipschitz on a neighborhood of $\left(0,z_0,(R_{2,0}(\lambda_{2,j})+\lambda_{2,j}A)_{1\le j\le N}\right)$.
For sufficiently small $h$, we extend the control by setting $u_t\mathrel{:=}0$ for $t\in[h,2h]$. The bounded and Lipschitz coefficients then give
\[
\sup_{t\in[h,2h]}
\left(
|\sigma_t^u-h^2|+|z_t^u-\widetilde z_h|
+|\mathbf R_{2,t}^u-\widetilde{\mathbf R}_{2,h}|
\right)=O(h).
\]
Together with $\widetilde z_h=z_0+O(\sqrt h)$ and
the expansion~\eqref{eq:negativity-variance-shift}, this gives
\begin{align}
\sigma_{2h}^u
&=h^2+\int_h^{2h}
\mu_\sigma(\sigma_t^u,z_t^u,\mathbf R_{2,t}^u)\,dt\notag\\
&=h^2+h\,\mu_\sigma(h^2,\widetilde z_h,\widetilde{\mathbf R}_{2,h})+O(h^2)\notag\\
&=h\,\mu_\sigma\!\left(
0,z_0,(R_{2,0}(\lambda_{2,j})+\lambda_{2,j}A)_{1\le j\le N}
\right)+O(h^{\frac{3}{2}}).
\label{eq:negativity-crossing-expansion}
\end{align}
For a fixed sufficiently small $h>0$ with $2h<T$,
the inequality~\eqref{eq:negativity-target-drift} and
the expansion~\eqref{eq:negativity-crossing-expansion} give $c:=\sigma_{2h}^u<0$.
Since $\beta_1<0$, $\Lambda_1>0$, and $c<0$, the denominator $\beta_1\Lambda_1c$ is nonzero. We extend the control to $[2h,T]$ by setting
\[
\forall t\in[2h,T],\qquad
u_t\mathrel{:=}-\frac{\mu_\sigma(c,z_t^u,\mathbf R_{2,t}^u)}{\beta_1\Lambda_1c}.
\]
This choice gives $\dot\sigma_t^u=0$, so $\sigma_t^u=c$, and
\begin{align*}
\forall t\in[2h,T],\ \forall j\in\{1,\ldots,N\},\qquad
R_{2,t}^u(\lambda_{2,j})
&=e^{-\lambda_{2,j}(t-2h)}R_{2,2h}^u(\lambda_{2,j})
+\left(1-e^{-\lambda_{2,j}(t-2h)}\right)c^2.
\end{align*}
This formula gives, for every $j\in\{1,\ldots,N\}$,
\[
\sup_{t\in[0,T]}R_{2,t}^u(\lambda_{2,j})
\le\max\left\{\sup_{t\in[0,2h]}R_{2,t}^u(\lambda_{2,j}),c^2\right\}<\infty.
\]
We write $\|\cdot\|_\infty$ for the supremum norm on $[0,T]$.
Since $\mathbf R_2^u$ is bounded and $\sigma^u=c$ on $[2h,T]$,
the equation~\eqref{eq:negativity-z-dynamics} and Grönwall's lemma give
\[
\|z^u\|_\infty+\|\mathbf R_2^u\|_\infty<\infty.
\]
The variance equations~\eqref{eq:negativity-variance-dynamics} and
the nonnegativity of $(\sigma^u)^2$ also give
\[
\forall t\in[0,T],\ \forall j\in\{1,\ldots,N\},\qquad
R_{2,t}^u(\lambda_{2,j})\ge e^{-\lambda_{2,j}t}R_{2,0}(\lambda_{2,j}).
\]
Summing with weights $\theta_{2,j}$ and using $\bar R_{2,0}>0$ from~\eqref{eq:initial-values} yields
\begin{equation}\label{eq:negativity-controlled-variance-lower-bound}
\inf_{t\in[0,T]}\bar R_{2,t}^u
\ge e^{-\lambda_{2,1}T}\bar R_{2,0}>0.
\end{equation}
The bound~\eqref{eq:negativity-controlled-variance-lower-bound} keeps the controlled path in the domain of $\mu_\sigma$. Together with $\beta_1\Lambda_1c\ne0$ and the boundedness of $z^u$ and $\mathbf R_2^u$, it ensures that the control on $[2h,T]$ is well defined and bounded. Thus $u\in L^\infty([0,T])\subset L^2([0,T])$
and $\sigma_T^u=c<0$.

\medskip

\noindent\textit{Step 2: Stability and the Cameron--Martin theorem.}
We set
\[
\forall t\in[0,T],\qquad
w_t\mathrel{:=}\int_0^t\left(u_s+\frac{\beta_1\Lambda_1}{2}\right)ds,
\qquad
G_t\mathrel{:=}\exp\left(\beta_1\Lambda_1(W_t-w_t)\right).
\]
The definition of $w$ gives
\[
G_t=\exp\left(\beta_1\Lambda_1W_t-\beta_1\Lambda_1\int_0^t u_s\,ds-\frac{\beta_1^2\Lambda_1^2}{2}t\right).
\]
Since $\dot w_t=u_t+\frac{\beta_1\Lambda_1}{2}$, It\^o's formula and \eqref{eq:negativity-sigma-dynamics} give
\begin{align*}
d(G_t^{-1})
&=(\beta_1\Lambda_1u_t+\beta_1^2\Lambda_1^2)
G_t^{-1}\,dt-\beta_1\Lambda_1G_t^{-1}\,dW_t,\\
d[\sigma,G^{-1}]_t
&=-\beta_1^2\Lambda_1^2\sigma_tG_t^{-1}\,dt.
\end{align*}
Since $G>0$, the ratio $\frac{\sigma}{G}$ is well defined and satisfies, for almost every $t\in[0,T]$,
\begin{equation}\label{eq:negativity-transformed-sigma}
\frac{d}{dt}\left(\frac{\sigma_t}{G_t}\right)
=G_t^{-1}\mu_\sigma(\sigma_t,z_t,\mathbf R_{2,t})
+\beta_1\Lambda_1u_t\frac{\sigma_t}{G_t}.
\end{equation}
For $G=1$, \eqref{eq:negativity-transformed-sigma} reduces to
\eqref{eq:negativity-sigma-dynamics} with $dW_t$ replaced by $u_t\,dt$.
By continuity of the controlled path and the lower bound~\eqref{eq:negativity-controlled-variance-lower-bound}, for sufficiently
small $\rho>0$, the open set
\[
\mathcal U_\rho:=\left\{(\sigma,z,\mathbf r):
\inf_{s\in[0,T]}\left(
|\sigma-\sigma_s^u|+|z-z_s^u|+|\mathbf r-\mathbf R_{2,s}^u|
\right)<\rho\right\}
\]
has compact closure in
$\{(\sigma,z,\mathbf r):\sum_{j=1}^N\theta_{2,j}r_j>0\}$.
The coefficients in \eqref{eq:negativity-lift} are bounded and Lipschitz
on $\overline{\mathcal U_\rho}$ by \eqref{eq:negativity-drift-definition}.
For $0<\varepsilon\le1$, on $\{\|W-w\|_\infty<\varepsilon\}$, we have
\begin{equation}\label{eq:negativity-G-stability}
\|G-1\|_\infty+\|G^{-1}-1\|_\infty\le C\varepsilon.
\end{equation}
Up to the first exit of $(\sigma,z,\mathbf R_2)$ from $\mathcal U_\rho$,
boundedness and Lipschitz continuity of $\mu_\sigma$, together with
the bound~\eqref{eq:negativity-G-stability}, give on $\overline{\mathcal U_\rho}$
\begin{align}
\left|
G_t^{-1}\mu_\sigma(\sigma_t,z_t,\mathbf R_{2,t})
-\mu_\sigma(\sigma_t^u,z_t^u,\mathbf R_{2,t}^u)
\right|
&\le
|G_t^{-1}-1|\,|\mu_\sigma(\sigma_t,z_t,\mathbf R_{2,t})|\notag\\
&\quad+C\left(
|\sigma_t-\sigma_t^u|+|z_t-z_t^u|
+|\mathbf R_{2,t}-\mathbf R_{2,t}^u|\right)\notag\\
&\le C\varepsilon+C\left(
|\sigma_t-\sigma_t^u|+|z_t-z_t^u|
+|\mathbf R_{2,t}-\mathbf R_{2,t}^u|\right).\label{eq:negativity-drift-comparison}
\end{align}
Moreover, the bound~\eqref{eq:negativity-G-stability} and boundedness of
$\sigma^u$ imply
\begin{align}
|\sigma_t-\sigma_t^u|
&=\left|
G_t\left(\frac{\sigma_t}{G_t}-\sigma_t^u\right)
+(G_t-1)\sigma_t^u
\right|\le C\left|\frac{\sigma_t}{G_t}-\sigma_t^u\right|+C\varepsilon.
\label{eq:negativity-sigma-comparison}
\end{align}
We compare \eqref{eq:negativity-transformed-sigma},
\eqref{eq:negativity-z-dynamics}, and
\eqref{eq:negativity-variance-dynamics} with the controlled equations.
Estimates~\eqref{eq:negativity-drift-comparison} and~\eqref{eq:negativity-sigma-comparison} give
\begin{align*}
\left|\frac{\sigma_t}{G_t}-\sigma_t^u\right|
+|z_t-z_t^u|
+|\mathbf R_{2,t}-\mathbf R_{2,t}^u|
&\le C\varepsilon\\
&\;{}+C\int_0^t(1+|u_s|)
\left(
\left|\frac{\sigma_s}{G_s}-\sigma_s^u\right|
+|z_s-z_s^u|
+|\mathbf R_{2,s}-\mathbf R_{2,s}^u|
\right)ds.
\end{align*}
Since $u\in L^2([0,T])\subset L^1([0,T])$, Grönwall's lemma
and the bound~\eqref{eq:negativity-sigma-comparison} give
\[
\sup_{s\in[0,t]}
\left(
|\sigma_s-\sigma_s^u|+|z_s-z_s^u|
+|\mathbf R_{2,s}-\mathbf R_{2,s}^u|
\right)\le C\varepsilon,
\]
where $C$ does not depend on $\varepsilon$.
For $C\varepsilon<\min\{\rho,-\frac{c}{2}\}$, the uniform bound on
$\sigma-\sigma^u$, $z-z^u$, and $\mathbf R_2-\mathbf R_2^u$,
together with continuity of $\sigma,z,\mathbf R_2$, excludes an exit
from $\mathcal U_\rho$ on $[0,T]$.
Hence the terminal volatility satisfies
\[
\sigma_T\le \sigma_T^u+C\varepsilon=c+C\varepsilon<0.
\]
Therefore the events satisfy
$\{\|W-w\|_\infty<\varepsilon\}\subset\{\sigma_T<0\}$.
Since
\[
\dot w=u+\frac{\beta_1\Lambda_1}{2}\in L^2([0,T]),
\]
the Cameron--Martin theorem and
$\mathbb P(\|W\|_\infty<\varepsilon)>0$ give
\[
\mathbb P(\sigma_T<0)
\ge
\mathbb P(\|W-w\|_\infty<\varepsilon)>0.
\]
The result in \cite[Corollary~2.6~(ii)]{GuyonJeanninJourdain2026Martingale}
applies with $t_0=\frac{T}{2}$, since
$\mathbb P(\sigma_{\frac{T}{2}}<0)>0$.
Hence $S$ is a strict local martingale on $[0,T]$
under $\mathbb P$.
\end{proof}

\bibliographystyle{plainnat}
\bibliography{bibtem}

\end{document}